\documentclass[12pt,a4paper]{article}

\usepackage{titlesec,titletoc}

\usepackage[margin=0.75in,width=6.25in]{geometry}
\usepackage{url}
\usepackage{fancyhdr}
\usepackage{setspace} %allows \onehalfspace spacing
\usepackage{appendix} %allows section appendices
\usepackage{enumerate}%allows custom enumeration values
\usepackage{calc}
\usepackage{xspace}
\usepackage{caption}
\usepackage{subcaption}
\usepackage{placeins}  %allows \FloatBarrier
\usepackage{pdflscape} % or {lscape}
\usepackage{longtable}
\usepackage{adjustbox} %fit a table into one page
\usepackage{xcolor,colortbl,booktabs}
\usepackage{parskip}  %no indent, parskip
\usepackage{authblk}  %more detailed author information
\usepackage{bold-extra} %allows bold small caps
\usepackage{stmaryrd} %use \llbracket 
\makeatletter
\def\thm@space@setup{%
\thm@preskip=\parskip \thm@postskip=0pt
}
\makeatother

\usepackage{amsmath}
\usepackage{amssymb}
\usepackage{amsthm}
\usepackage{mathtools}
\usepackage{float}

\usepackage{algorithm}
\usepackage{algorithmic}
\usepackage{tikz}
\usetikzlibrary{decorations.pathreplacing}
\usepackage{pgfkeys}

\usepackage{graphicx}
\graphicspath{ {./figures/} }

\usepackage{listings}
\usepackage{xcolor}

\lstdefinestyle{mystyle}{
	backgroundcolor=\color{lightgray}, % Set background color
	commentstyle=\color{green},
	keywordstyle=\color{blue},
	numberstyle=\tiny\color{gray},
	stringstyle=\color{red},
	basicstyle=\ttfamily\footnotesize,
	breakatwhitespace=false,
	breaklines=true,
	captionpos=b,
	keepspaces=true,
	numbers=left,
	numbersep=5pt,
	showspaces=false,
	showstringspaces=false,
	showtabs=false,
	tabsize=2
}

\newcommand{\mcf}{\mathcal}

\newcommand{\tpose}{^{T}}

\renewcommand{\Re}{\textbf{R}}
\newcommand{\reals}{\Re}
\renewcommand{\SS}{\textbf{S}}

\newcommand{\eqdef}{=}

\newcommand{\norm}[1]{\left\| #1 \right\|}
\newcommand{\set}[2]{\left\{ #1\ \left| \ #2 \right. \right\}}

\newcommand{\innerprod}[2]{\left\langle{#1},{#2}\right\rangle}

\theoremstyle{remark}

\usepackage{xcolor,calc}

\newcommand{\detailtablecaption}{FOO}

\definecolor{BenchHighlight}{rgb}{0.8,0.8,0.9}
\newcommand{\winner}{\cellcolor{BenchHighlight}}

\title{\LARGE  {\bfseries CuClarabel}: \textbf{GPU Acceleration for a Conic Optimization Solver}}
\author[1]{Yuwen Chen}  
\author[2]{Danny Tse}
\author[3]{Parth Nobel}
\author[1]{Paul Goulart}
\author[3]{Stephen Boyd}

\affil[1]{Department of Engineering Science, University of Oxford,
	Oxford, UK}
\affil[2]{Department of Computer Science, Stanford University, Stanford, CA}
\affil[3]{Department of Electrical Engineering, Stanford University, Stanford, CA}
\date{\today}
\begin{document}
%\keywords{GPU computing, interior point method, conic optimization}

\maketitle
 \begin{abstract}
We present the GPU implementation of the general-purpose interior-point solver Clarabel
for convex optimization problems with conic constraints.
We introduce a mixed parallel computing strategy that processes
linear constraints first, then handles other conic constraints in parallel.
{
The GPU solver currently supports linear equality and inequality constraints, second-order cones, exponential cones, power cones and positive semidefinite cones of the same dimensionality.}
We demonstrate that integrating a mixed parallel computing strategy with 
GPU-based direct linear system solvers enhances the performance of GPU-based conic solvers, 
surpassing their CPU-based counterparts across a wide range of conic optimization problems. 
We also show that employing mixed-precision linear system solvers can potentially achieve 
additional acceleration without compromising solution accuracy.
\end{abstract}

%------------------------------
% Paper body
%------------------------------

%\tableofcontents
\section{Introduction}
We consider the following convex optimization problem~\cite{bv_2004} with a quadratic objective and conic constraints:
\begin{equation}\label{eqn:primal}\tag{$\mcf P$}
\begin{array}{ll}
\text{minimize}   & \frac{1}{2}x\tpose P x + q\tpose x \\
\text{subject to} & Ax + s = b, \\
		  & s \in \mcf{K},
\end{array}
\end{equation}
with respect to $x, s$ and
with parameters
$A\in \Re^{m\times n}$, $b\in \Re^m$, $q \in \Re^n$ and $P \in \SS_+^n$ and variables $x\in \Re^n, s\in \Re^m$.
For the rest of the paper, $\SS_+^n$ represents the cone of positive semi-definite matrices.
The cone $\mcf{K}$ is a closed convex cone.
The formulation~\eqref{eqn:primal} is very general and can model most conic convex optimization problems in practice.
Examples include the optimal power flow problem in power systems~\cite{PGLIB:2019}, model predictive control in 
control~\cite{Borrelli:MPCbook, lmi_boyd}, limit analysis of engineering structures in mechanics~\cite{Makrodimopoulos2007}, 
support vector machines \cite{Cortes:SVM:1995} and lasso problems~\cite{Tibshirani:1996} in machine learning, statistics, 
and signal processing, and portfolio optimization in finance~\cite{Krokhmal2007,Vinel2017}.
Linear equality (zero cones) and inequality (nonnegative cones), second-order cone~\cite{socp_boyd}, and semidefinite cone~\cite{sdp_boyd} constraints
have been long supported in standard conic optimization solvers, and support for exponential and power cone constraints was recently included in several 
state-of-the-art conic optimization solvers~\cite{mosek,scs,COSMO,Clarabel}.
The combination of these cones can represent many more elaborate convex constraints through the lens of disciplined convex programming~\cite{Grant04, cvx, cvxr}.

The dual problem of~\eqref{eqn:primal} is
\begin{equation}\label{eqn:dual}\tag{$\mcf{D}$}
\begin{array}{ll}
\text{maximize}    & -\frac{1}{2}x\tpose P x - b\tpose z \\
\text{subject to}  & P x  +  A\tpose z = -q, \\
                   & z \in \mcf{K}^*,
\end{array}
\end{equation}
with respect to $x, z$ and where $\mathcal K^*$ is the dual cone of $\mathcal K$. 
Solving~\ref{eqn:primal} and~\ref{eqn:dual} is equivalent to solving 
the Karush-Kuhn-Tucker (KKT) conditions when strong duality holds. On the other hand,
the set of strongly primal infeasibility certificates for \eqref{eqn:primal} is
\begin{align}
	\mathbb{P} &\eqdef \set{z}{ A\tpose z = 0,~z\in \mcf{K}^*,~\innerprod{b}{z} < 0},
\end{align}
and the set of strongly dual infeasibility certificates is
\begin{align}
	\mathbb{D} &\eqdef \set{x}{~Px = 0,~-Ax \in \mcf{K},~\innerprod{q}{x} < 0}. 
\end{align} 
Finding an optimal solution with strong duality or detecting a strongly infeasible certificate 
can be unified into solving a linear complementarity problem 
with two additional slack variables $\tau, \kappa \ge 0$~\cite{scsqp}, 
which can be reformulated to the following problem:
\begin{equation}\label{eqn:PD_homogeneous}\tag{$\mcf H$}
\begin{array}{ll}
\text{minimize}   & s^Tz+ \tau\kappa \\
\text{subject to} & \frac{1}{\tau}x\tpose P x + q\tpose x + b\tpose z = -\kappa, \\
		  & Px + A\tpose z + q\tau = 0, \\
		  & Ax + s - b\tau = 0, \\
		  & (s,z,\tau,\kappa) \in \mcf{K} \times \mcf{K^*} \times \Re_+ \times \Re_+,
\end{array}
\end{equation}
with respect to $x, s, z, \tau, \kappa$. {It is well-known that when $P = 0$, (\ref{eqn:PD_homogeneous}) is 
the \emph{homogenous self-dual embedding}~\cite{hsde} of (\ref{eqn:primal}), 
and the extension for $P \ne 0$ can be regarded as a \emph{homogeneous embedding}
for linear complementarity problems~\cite{Andersen1999,Yoshise07}. Precisely, in our case}
 (\ref{eqn:PD_homogeneous}) is homogeneous but not self-dual.
In~\cite{Clarabel} it is also shown that \ref{eqn:PD_homogeneous} is always (asymptotically) feasible,
and we can recover either an optimal solution or a strong infeasibility certificate 
of $\eqref{eqn:primal}$ and $(\mcf{D})$,
depending on the value of the optimal solution 
$(x^\star,z^\star,s^\star,\tau^\star,\kappa^\star)$ to \eqref{eqn:PD_homogeneous}:
\begin{enumerate}[{i)}]
	\item If $\tau^\star > 0$ then $(x^\star/\tau^\star,s^\star/\tau^\star)$
	is an optimal solution to \eqref{eqn:primal} and
	$(x^\star/\tau^\star,z^\star/\tau^\star)$ is an optimal solution to $(\mcf{D})$.
	\item If $\kappa^\star > 0$ then at least one of the following holds:
	\begin{itemize}
		\item \eqref{eqn:primal} is strongly infeasible and $z^\star \in \mathbb{P}$.
		\item $(\mcf{D})$ is strongly infeasible and $x^\star \in \mathbb{D}$.
	\end{itemize}
\end{enumerate}
The optimal solution $\tau^\star, \kappa^\star$ satisfy the complementarity slackness condition,
\emph{i.e.}, at most one of $\tau^\star, \kappa^\star$ is nonzero. 
The pathological case
$\tau^\star=\kappa^\star=0$ 
has been discussed in~\cite{scsqp}.

The \emph{interior-point method}~\cite{nesterov_nemirovskii} is a popular choice for solving $\eqref{eqn:PD_homogeneous}$.
However, it usually requires to factorize linear systems that are increasingly ill-conditioned.
{The time complexity of matrix factorization scales with respect to the fill-in that is closely related to nonzeros of the matrix and permutation strategies during factorization.} 
It is thus time-consuming to solve large-scale conic optimization problems with interior-point methods.

\subsection{Our contribution}\label{subsec:contribution}
We contribute a GPU-accelerated implementation of the general-purpose interior-point solver Clarabel~\cite{Clarabel}, in Julia.
This implementation  supports cases where $\mcf{K}$ is a product of
\emph{atomic} cones: 
{
zero cones, nonnegative cones, second-order cones, exponential cones, power cones and positive semidefinite cones. 
}
We propose a \emph{mixed parallel computing strategy} that parallelizes computing for each type of cone,
integrates the cuDSS~\cite{CUDSS.jl} library for linear system solving,
and supports \emph{mixed-precision} linear system solves for moderate speed improvements.
Furthermore, we evaluate our solver against others across a variety of conic optimization problems. {Our implementation of CuClarabel is open-sourced on GitHub\footnote{\url{https://github.com/oxfordcontrol/Clarabel.jl/tree/CuClarabel}} and can be easily accessed through CVXPY~\cite{cvxpy}.}

\subsection{Related work}\label{subsec:related-work}

\paragraph*{Interior-point methods and solver development.}
The interior-point method~\cite{nesterov_nemirovskii} was first discovered by Dikin~\cite{dikin}, and became more mainstream after the conception of Karmarkar's method~\cite{karmarkar}, 
a polynomial-time algorithm for linear programming, and Renegar's~\cite{renegar} path following method. Interior-point methods are known for their ability to solve conic optimization problems to high precision,
and is chosen as the default algorithm for many conic optimization solvers~\cite{mosek,ecos,Clarabel,Hypatia,DDS}. Common variations of the interior-point method include potential reduction methods and path-following methods.

Interior-point methods employ a Newton-like strategy to compute a search direction at every iteration.   Unfortunately, the matrix factorization required 
to compute this direction which scales with the dimension of an optimization problem, rendering very large problems difficult to solve.
Developing efficient interior-point methods on exotic cones directly is a promising research direction~\cite{Hypatia,DDS} to alleviate this computational burden.
Using exotic cones, we can represent equivalent problems with significantly fewer variables and exploit sparse structure within these exotic cones 
for efficient implementation of interior-point methods~\cite{ecos,Andersen10,Chen2023}. However, operations within an exotic cone can hardly be parallelized, 
while parallelism across heterogeneous cones of different dimensionalities will introduce significant synchronization delay.

\paragraph*{GPU acceleration in optimization algorithms.}
GPUs are playing an increasingly significant role in scientific computing.
In optimization, GPUs are used to run solvers based on first-order methods.
For example, CuPDLP~\cite{cupdlp} is based on the popular PDHG algorithm~\cite{pdhg},
which requires only matrix multiplication and addition without the use of direct methods (\emph{i.e.},\ it is \emph{factorization-free}).
Solvers that require the solution to linear systems, like SCS~\cite{scs} and CuOSQP~\cite{cuosqp},
have relied on indirect iterative methods---such as the conjugate gradient (CG),
the minimal residual (MINRES)~\cite{Saad2003},
and the generalized minimum residual (GMRES)~\cite{gmres} methods---to solve these systems on GPUs.
However, the linear systems solved in first-order methods are generally much better conditioned than those encountered in interior-point methods, 
where the linear systems become increasingly ill-conditioned as the iterations progress. As an interior-point method approaches higher precision, 
the number of iterations for each inner indirect linear solves increases significantly, which will eventually offset benefits of GPU parallelism 
and make GPU-based interior-point methods less preferable compared to CPU-based solvers with direct methods in overall computational time~\cite{Smith2012}.	
Recently, NVIDIA released the cuDSS package~\cite{CUDSS}, which provides fast direct methods on GPUs for sparse linear systems.
Previous work has integrated cuDSS in a nonlinear optimization solver, MADNLP~\cite{Pacaud2024}, resulting in significant speed-up on large-scale problems 
compared to its CPU-based counterpart. 

\paragraph*{Mixed-precision methods.}
Mixed-precision, or multiprecision, methods~\cite{survey_mixed_prec, Higham_Mary_2022} have become progressively more popular due to their synergies with modern GPU architectures.
For example, in lower-precision configuration one enjoys significant speedup in algorithms, as modern architectures have 32-bit implementations that are around twice as fast as their 64-bit counterparts~\cite{survey_mixed_prec, accelerating_scientific_computations}.
Mixed precision methods aim to capitalize on the computational benefits of performing expensive operations in lower precision, while maintaining (or reducing the negative impact to) the superior numerical accuracy from higher precision.
Classically, mixed-precision methods for direct (linear system solving) methods involve factorizing a matrix in lower precision, and then applying iterative refinement.
In this approach, only the more expensive steps, such as the matrix factorizations and backsolves, are done in lower precision.
Mixed-precision methods are used throughout scientific computing, with applications including BLAS operations~\cite{blas_mixed}, different linear system solvers such as Krylov methods (CG, GMRES)~\cite{Higham_Mary_2022, carson2021mixed, gratton2020exploiting}, solving partial differential equations~\cite{hayford2024speeding}, and training deep neural networks~\cite{hayford2024speeding, micikevicius2018mixed, courbariaux2015training}.

\subsection{Paper outline}

In \S\ref{sec:clarabel-review}, we review supported cones and discuss the interior point method used by Clarabel~\cite{Clarabel} for solving conic optimization problems.
In \S\ref{sec:GPU-cone}, we outline how to implement parallel computation for cone operations and solve linear systems within our GPU solver.
\S\ref{sec:numerical} details our numerical experiments.
We detail the scaling matrices for each cone in~\S\ref{appendix:scaling-matrix}.

\section{Backgrounds}\label{sec:clarabel-review}
{
We first review supported cones in CuClarabel with the corresponding barrier functions, along with their dual cones.
We then briefly sketch the main operations used by the Clarabel solver~\cite{Clarabel} to compute a solution to \eqref{eqn:primal}.  
}

\subsection{Supported cones}
We support the following atomic cones in our GPU solver:
\begin{itemize}
	\item
	The \emph{zero cone}, defined as 
	$$
	\{0\}^n\eqdef\left\{x \in \Re^n \ \middle| \ x_i=0,~ i=1, \ldots, n\right\}.
	$$
	The dual cone of the zero cone is $\left(\{0\}^n\right)^*=\Re^n$.
	
	\item
	The \emph{nonnegative cone}, defined as
	$$
	\Re_{+}^n\eqdef\left\{x \in \Re^n \ \middle| \ x_i \geq 0, ~ i=1, \ldots, n\right\}.
	$$
	The nonnegative cone is a self-dual convex cone, \emph{i.e.}, $\left(\Re_{+}^{n}\right)^*=\Re_{+}^n$.
	\item 
	The \emph{second-order cone} $\mathcal{K}_{\text {soc }}^n$ (also called the \emph{quadratic} or \emph{Lorentz cone}), defined as
	$$
	\mathcal{K}_{\mathrm{soc}}^n\eqdef\left\{(t, x) \ \middle| \ x \in \Re^{n-1},~ t \in \Re_{+},~\|x\|_2 \leq t\right\}.
	$$
	The second-order cone is self-dual, \emph{i.e.}, $\mathcal{K}_{\mathrm{soc}}=\mathcal{K}_{\mathrm{soc}}^*$.
	\item 
{
	The \emph{positive semidefinite cone} $\mathcal{K}_{\mathrm{psd}}^n$ is defined as
	$$
	\mathcal{K}_{\mathrm{psd}}^n \eqdef \left\{ X \in \mathbb{S}^n \ \middle| \ X \succeq 0 \right\},
	$$
	where $\mathbb{S}^n$ denotes the space of $n \times n$ symmetric matrices, and $X \succeq 0$ indicates that $X$ is positive semidefinite. The positive semidefinite cone is self-dual, \emph{i.e.}, $\mathcal{K}_{\mathrm{psd}}^n = \left(\mathcal{K}_{\mathrm{psd}}^n\right)^*$.\footnote{{In the current solver implementation, all positive semidefinite (PSD) cones are restricted to have dimension at most $32$, and all PSD cones within a given problem must share the same dimension.}}
}

	\item
	The \emph{exponential cone}, a $3$-dimensional cone defined as
	$$
	\mathcal{K}_{\exp }\eqdef\left\{(x, y, z) \ \middle| \ y>0,~ y 
	\exp \left(\frac{x}{y}\right) \leq z\right\} \cup\{(x, 0, z) \mid x \leq 0, ~z \geq 0\},
	$$
	with its dual cone given by
	$$
	\mathcal{K}_{\exp }^*=\left\{(u, v, w) \ \middle| \ u<0,~-u \exp \left(\frac{v}{u}\right) \leq \exp (1) w\right\} \cup\{(0, v, w) \mid v \geq 0,~ w \geq 0\}.
	$$
	
	\item
	The $3$-dimensional \emph{power cone} with exponent $\alpha \in (0,1)$, defined as
	$$
	\mathcal{K}_{\text{pow},\alpha} = \left\{(x,y,z) \ \middle| \ x^\alpha y^{1-\alpha} \geq \vert z \vert, ~x \geq 0,~ y \geq 0\right\},
	$$
	with its dual cone given by
	$$
	\mathcal{K}_{\text{pow}, \alpha}^*=\left\{(u, v, w) \ \middle| \ \left(\frac{u}{\alpha}\right)^\alpha\left(\frac{v}{1-\alpha}\right)^{1-\alpha} \geq \vert w \vert,~ u \geq 0,~ v \geq 0\right\}.
	$$
\end{itemize}

\subsection{Interior point method}
Solving the problem \eqref{eqn:PD_homogeneous} amounts to finding a root of the following nonlinear equations 
\begin{equation}\label{eqn:root_finding_G}
\begin{aligned}
		G(x,z,s,\tau,\kappa) \eqdef 
	\begin{bmatrix}
		0\\s\\ \kappa 
	\end{bmatrix}
	- &
	\begin{bmatrix}
		\hphantom{+}P  & A\tpose & q \\
		-A & 0 & b\\
		-q\tpose & -b\tpose & 0
	\end{bmatrix}
	\begin{bmatrix}
		x \\ z \\ \tau 
	\end{bmatrix}
	+ 
	\begin{bmatrix}
		0 \\ 0 \\
		\frac{1}{\tau}x\tpose P x
	\end{bmatrix} = 0, \\
	& (x,z,s,\tau,\kappa) \in \mathcal{F},
\end{aligned}
\end{equation}
{where $\mathcal{F} : \eqdef  \Re^n \times \mathcal{K}^* \times \mathcal{K} \times\Re_+ \times \Re_+$ 
	defines the region of cone constraints.
} 
Besides the zero cone that is a linear constraint, 
other supported conic constraints are smoothed by nonlinear equations in pairs within an interior-point method,
\begin{align}
	s=-\mu \nabla f(z), \quad \tau \kappa = \mu,
\end{align} 
where $\mu$ is the smoothing parameter and $f(\cdot)$ is the logarithmically homogeneous self-concordant barrier (LHSCB) function for cone $\mathcal{K}^*$.
The barrier functions for different cones are:
{
\begin{enumerate}
	\item \textit{Nonnegative cone} $\reals_{+}^n$ of degree $n$:
	\begin{equation}
		f(z) = -\sum_{i \in \llbracket n \rrbracket}\log(z_i), \ z \in \reals_{+}^n.
	\end{equation}
	\item \textit{Second order cone} $\mathcal{K}_{q}^n$ of degree $1$:
	\begin{equation}
		f(z) = -\frac{1}{2}\log\left(z_1^2 - \sum_{i =2}^n z_i^2\right), \ z \in \mathcal{K}_{q}^n. \label{eq:soc-barrier}
	\end{equation}
	{
	\item \textit{Positive semidefinite cone} $\mathcal{K}_{\mathrm{psd}}^n$ of degree $n$:
	\begin{equation}
		f(X) = -\log\det(X), \quad X \in \mathcal{K}_{\mathrm{psd}}^n. \label{eq:psd-barrier}
	\end{equation}
	}
	\item \textit{Dual exponential cone} $\mathcal{K}_{\mathrm{exp}}^*$ of degree $3$:
	\begin{equation}
		f(z) = -\log\left(z_2 - z_1 - z_1\log\left(\frac{z_3}{-z_1}\right)\right) - \log(-z_1) - \log(z_3), \ z \in \mathcal{K}_{\mathrm{exp}}^*.
	\end{equation}
	\item \textit{Dual power cone} $\mathcal{K}_{\mathrm{pow}}^*$ of degree $3$:\\
	\begin{small}
		\begin{equation}
			f(z) = -\log\left(\left(\frac{z_1}{\alpha}\right)^{2\alpha} \left(\frac{z_2}{1-\alpha}\right)^{2(1-\alpha)} - z_3^2\right) - (1-\alpha)\log(z_1) - \alpha\log(z_2), \ z \in \mathcal{K}_{\mathrm{pow}}^*.
		\end{equation}
	\end{small}
\end{enumerate}
}

The trajectory (also called the \textit{central path})
\begin{align}
\begin{aligned}
		G(v) = \mu G(v^0), \\
	s=-\mu \nabla f(z), \quad \tau \kappa = \mu,
\end{aligned}\label{eqn:central-path}
\end{align}
where $v \eqdef(x,z,s,\tau,\kappa)$,
characterizes the solution of~\eqref{eqn:root_finding_G} in the right limit $\mu \rightarrow 0^+$, given an initial point $v^0$. 
After starting from $v^0$, the Clarabel solver iterates the following steps for each iteration $k$:

\paragraph*{Update residuals and check the termination condition.} 
We update several key metrics at the start of each iteration.
Defining the normalized variables $\bar{x} = x/\tau, \bar{s} = s/\tau, \bar{z} = z/\tau$, the primal and dual residuals are then
\begin{align}
	\begin{aligned}
		r_p &\eqdef -A\bar x - \bar s + b, \\
		r_d &\eqdef P\bar x + A^T \bar z + q, \\
		g_p &\eqdef \frac{1}{2} \bar x^T P \bar x + q^T \bar x, \\
		g_d &\eqdef -\frac{1}{2} \bar x^T P \bar x - b^T \bar z
	\end{aligned}\label{eq:residuals-objectives}
\end{align}

and complementarity slackness $\mu = \frac{{s}^T z +\kappa \tau}{\nu + 1}$, where $\nu$ is the degree of cone $\mathcal{K}$.
The solver returns an approximate optimal point if
\begin{align}
\begin{aligned}
	\norm{r_p}_\infty &< \epsilon_{\text{feas}} \max\{1, \norm{b}_\infty + \norm{\bar x}_\infty + \norm{\bar s}_\infty\} \\
	\norm{r_d}_\infty &< \epsilon_{\text{feas}} \max\{1, \norm{q}_\infty + \norm{\bar x}_\infty + \norm{\bar z}_\infty\} \\
	|g_p - g_d| &< \epsilon_{\text{feas}} \max\{1, \min\{|g_p|, |g_d|\}\}.
\end{aligned}\label{eq:termination-check}
\end{align}
Otherwise, it returns a certificate of primal infeasibility if
\begin{subequations}\label{eq:infeasibility-check}
\begin{align}
\begin{aligned}
	\norm{A\tpose z}_\infty & < -\epsilon_{\text{inf}} \max(1, \norm{x}_\infty + \norm{z}_\infty) (b\tpose z) \\
	b\tpose z & < -\epsilon_{\text{inf}}, \\
\end{aligned}
\end{align}
\\
\text{and a certificate of dual infeasibility if }
\\
\begin{align}
\begin{aligned}
\norm{Px}_\infty & < -\epsilon_{\text{inf}} \max(1, \norm{x}_\infty) (b\tpose z) \\
\norm{Ax + s}_\infty & < -\epsilon_{\text{inf}} \max(1, \norm{x}_\infty + \norm{s}_\infty) (q\tpose x) \\
q\tpose x  & < -\epsilon_{\text{inf}}.
\end{aligned}
\end{align}
\end{subequations}
Note that $\epsilon_{\text{feas}},\epsilon_{\text{inf}}$ are predefined parameters within the Clarabel solver.

\paragraph*{Find search directions.}
We then compute Newton-like search directions using a linearization of the central path~\eqref{eqn:central-path}.
In other words, we solve the following linear system given some right-hand side residual $d \eqdef (d_x, d_z,d_\tau,$ $d_s,d_\kappa)$,
\begin{subequations}\label{eqn:linsys_4x4}
	\begin{gather}%
		\begin{bmatrix}
			0 \\ \Delta s \\ \Delta \kappa
		\end{bmatrix}  
		-
		\begin{bmatrix}
			P & A\tpose & q \\
			-A  & 0 & b \\
			-(q + 2P\xi)\tpose & -b\tpose & \xi\tpose P \xi
		\end{bmatrix}
		\begin{bmatrix}
			\Delta x \\ \Delta z \\ \Delta \tau
		\end{bmatrix}  
		= 
		-
		\begin{bmatrix}
			d_x\\d_z \\ d_\tau
		\end{bmatrix}  \label{eqn:linsys_4x4:a} \\[1ex] 
		H \Delta z + \Delta s = -d_s, \quad \kappa \Delta \tau + \tau\Delta \kappa = - d_\kappa,
		\label{eqn:linsys_4x4:b}
	\end{gather}
\end{subequations}
where $\xi = x\tau^{-1}$, and $H$ is the same positive-definite scaling matrix as in the CPU version of Clarabel~\cite{Clarabel}, also described in \S\ref{appendix:scaling-matrix}.
We have shown in~\cite{Clarabel} that solving \eqref{eqn:linsys_4x4} reduces to solve the next linear system with two different right-hand sides,
\begin{equation}\label{eqn:linsys_2x2}
	\underbrace{\begin{bmatrix}
			P & A^T\\
			A & - H
	\end{bmatrix}}_{K}
	\left[
	\begin{array}{c|c}
		\Delta x_1 & \Delta x_2\\
		\Delta z_1 & \Delta z_2
	\end{array}
	\right]
	=
	\left[
	\begin{array}{c|c}
		d_x  & -q \\
		-(d_z - d_s) & b
	\end{array}
	\right].
\end{equation}
{Since both $P,H$ are positive semidefinite, adding a positive regularization to them makes $K$ symmetric quasi-definite and thus strongly factorizable via the $LDL^T$ decomposition~\cite{vanderbei_sqd}.}
After solving~\eqref{eqn:linsys_2x2}, we recover the search direction $\Delta \eqdef (\Delta x,\Delta z,\Delta \tau, \Delta s, \Delta \kappa)$ using 
\begin{subequations}\label{eqn:step-direction}
	\begin{gather}
		\Delta \tau =\frac{d_\tau-d_\kappa /\tau+(2P\xi + q)^T \Delta x_1+b^T \Delta z_1}{\kappa/\tau + \xi^T P \xi-(2P\xi + q)^T \Delta x_2-b^T \Delta z_2} \notag \\[1ex]
		=\frac{d_\tau-d_\kappa /\tau + q^T \Delta x_1+b^T \Delta z_1 + 2\xi^T P \Delta x_1}{\|\Delta x_2 - \xi\|_P^2 -\|\Delta x_2\|_P^2 -q^T \Delta x_2-b^T \Delta z_2},
		\intertext{and }
		 \Delta x=\Delta x_1+\Delta \tau \Delta x_2, \quad \Delta z=\Delta z_1+\Delta \tau \Delta  z_2,\\[1ex]
		\Delta s = -d_s - H \Delta z, \quad \Delta \kappa=-(d_\kappa+\kappa \Delta \tau) /\tau.
	\end{gather}
\end{subequations}
In an interior-point method with a predictor-corrector scheme, we need to solve~\eqref{eqn:linsys_4x4} with two different values for $d$.   The first is for the affine step (predictor) with
\begin{align}
	d = (G(x,z,s,\tau,\kappa), \kappa \tau, s). \label{eq:affine-step}
\end{align}
The other assigns $d$ as
\begin{align}
	\begin{aligned}
		(d_x, d_z, d_{\tau}) = (1-\sigma)G(x,z,s,\tau,\kappa), \quad d_{\kappa} = \kappa \tau + \Delta \kappa \Delta \tau - \sigma \mu, \\
		d_s = \left\{
		\begin{array}{cl}
			W^T \left(\lambda \backslash \left(\lambda \circ \lambda + \eta - \sigma \mu \mathbf{e}\right) \right) & \text{(symmetric)}\\
			s + \sigma \mu \nabla f(z) + \eta & \text{(nonsymmetric),}
		\end{array}
		\right.
	\end{aligned}\label{ds-correction-step}	
\end{align} 
in the combined step (predictor+corrector), where $\lambda \eqdef W^{-T}s = W z$ and $\mathbf{e}$ is the idempotent for a symmetric cone with the product operator `$\circ$' and its inverse operator `$\backslash$'~\cite{cvxopt}.
Here $\eta$ denotes a higher-order correction term, which is a heuristic technique that can significantly accelerate the convergence of interior-point methods~\cite{Wright97}.
We set it to the Mehrotra correction~\cite{Mehrotra1992} 
\begin{subequations}\label{eqn:higher-order-correction}
	\begin{gather}
	\eta = (W^{-1} \Delta s) \circ (W \Delta z),
	\intertext{for symmetric cones and the $3$rd-order correction~\cite{Dahl21}}
	\eta = -\frac{1}{2} \nabla^3 f(z)[\Delta z^a,\nabla^2 f(z)^{-1} \Delta s^a],
	\end{gather}
\end{subequations}
for nonsymmetric cones.
{
The centering parameter $\sigma$ controls the rate at which both the residual $G(x, z, s, \tau, \kappa)$ and the complementarity measure $\mu$ decrease. 
It is determined heuristically based on the value of the affine step size $\alpha_a$, defined as the maximal step size ensuring that $v + \alpha_a \Delta_a \in \mathcal{F}$. Computing $\alpha_a$ is non-trivial in the presence of exponential and power cones. As a practical alternative, we perform a backtracking line search with a shrinking ratio of $0.8$ to determine a feasible $\alpha_a$ for power and exponential cones.}

\paragraph*{Update iterates.}
At the end of each iteration $k$, we move the current iterate $v$ along the combined direction $\Delta_c := (x,z,s,\tau,\kappa)$ and obtain the new iterate $v + \alpha_c \Delta_c$.
{
As discussed, the combined step size $\alpha_c$ should satisfy $v + \alpha_c \Delta_c \in \mathcal{F}$.
Furthermore, we must ensure that the new iterate $v + \alpha_c \Delta_c$ 
stays in the neighborhood $\mathcal{N}(\beta)$ of the central path~\eqref{eqn:central-path}~\cite{Dahl21},
\begin{align}
	\mathcal{N}(\beta) = \left\{ (s,z) \in \mathcal{K} \times \mathcal{K}^* \;\middle|\; \nu_i \left\langle \nabla f^*(s_i), \nabla f(z_i) \right\rangle^{-1} \geq \beta \mu, \; i = 1, \ldots, p \right\},
\end{align}
where $\beta$ is set to $10^{-6}$.
}

\section{GPU formulation of Clarabel}\label{sec:GPU-cone}

\subsection{A primer on GPU programming}

Originally developed for computationally demanding gaming applications, graphics processing units (GPUs) are \emph{massively parallel}, \emph{multithreaded}, \emph{manycore processors}, 
with massive computational power and memory bandwidth. As a result, GPUs are used throughout scientific computing. In this section, we outline some core properties of 
GPU computing devices to highlight why a GPU implementation is a natural extension to Clarabel. We then discuss our specific implementation details.

Due to their natural parallelism, GPUs differ dramatically from CPUs in the way their transistors are configured. GPUs have a smaller number of caches (blocks where memory access is fast) 
and instruction processing blocks, and far more, albeit simpler, computational blocks (\emph{i.e.},\ arithmetic logic and floating point units). This suggests that CPUs utilize their larger caches 
to minimize instruction and memory latency within each thread, while GPUs switch between their significantly larger number of threads to hide said latency.

GPUs use a \emph{single instruction, many threads} (SIMT) approach. In practice, the GPU receives a stream of instructions -- each instruction is sent to groups of GPU cores 
(also known as a \emph{warp}) and acts on multiple data in parallel. Each group of GPU cores, as a result, has \emph{single instruction, many data} (SIMD) structure. 
This contrasts the traditional CPU vector lane approach of \emph{single instruction, single data}. Note that modern CPUs also support SIMD, but at a much smaller scale than GPUs, 
as implied by the transistor layout. Furthermore, the SIMD structure requires data parallelism for parallel execution on a GPU.

{
\subsection{GPU implementation for interior point method}
We detail our GPU implementation of the primal-dual interior point method in Algorithm~\ref{alg:IPM}, which is divided into two phases: 
setup and solve. 
The setup phase involves data equilibration and solver initialization, 
while the solve phase implements the classical primal-dual interior point method. 
Synchronization is required at the end of each step in Algorithm~\ref{alg:IPM}.}
\begin{algorithm}
	\caption{Primal-dual interior-point method}
	\begin{algorithmic}[1]
		\REQUIRE Parameter inputs $P,A,q,b$ and cone information $\mathcal{K}$. \\
		\vspace*{2mm}
		//Setup phase
		\STATE Equilibrate $P,A,q,b$.
		\STATE Initialize solver structure $s$, including matrix $K$ and cuDSS solver.\\
		\vspace*{2mm}
		//Solve phase
		\WHILE{termination check~\eqref{eq:termination-check} or~\eqref{eq:infeasibility-check} is not satisfied}
		\STATE Update the scaling matrix $H$ as in Appendix~\ref{appendix:scaling-matrix}.
		\STATE Factorize the matrix $K$ in~\eqref{eqn:linsys_2x2}.
		\STATE Compute the right-hand residual $d$~\eqref{eq:affine-step} for the affine step 
		and solve~\eqref{eqn:linsys_2x2} to obtain $\Delta_a$ via~\eqref{eqn:step-direction}.
		\STATE Compute the affine step $\alpha_a$ satisfying $v + \alpha_a \Delta_a \in \mathcal{F}$.
		\STATE Compute the centering parameter $\sigma = (1-\alpha_a)^3$.
		\STATE Compute the right-hand residual $d$~\eqref{ds-correction-step} for the combined step 
		and solve~\eqref{eqn:linsys_2x2} to obtain $\Delta_c$ via~\eqref{eqn:step-direction}.
		\STATE Compute the combined step $\alpha_c$ satisfying $v + \alpha_c \Delta_c \in \mathcal{F}$ with the neighborhood check.
		\STATE Update the variable $v \leftarrow v + 0.99 \alpha_c \Delta_c$.
		\STATE Update residuals and objective values in~\eqref{eq:residuals-objectives}.
		\ENDWHILE
	\end{algorithmic}\label{alg:IPM}
\end{algorithm}

{
The matrix factorization (step 5) and the back-solve (steps 6 and 9) are the most time-consuming parts in each iteration of an interior point method, and can be computed using the cuDSS package. Variable and information update (steps 11 and 12) contain matrix addition and multiplication operations that have already been supported in CUDA. Steps 4, 6, 7, 9, and 10 include cone operations that should be tailored for each cone, but we can abstract them into the same framework for parallelism. We will detail how to parallelize each of these steps under a unified framework in the next subsection.
}

\subsection{Mixed parallel computing strategy}\label{subsec:kernel-operations}
The cone operations related to steps 4, 6, 7, 9, and 10 above require only information local to each constituent cone, and hence can be executed concurrently 
with respect to individual cones. We propose the \emph{mixed parallel computing strategy} for a cone operation across different types of cones. 
Each family of cones is handled in parallel, and families of cones {can be addressed by separate streams}. 

As stated earlier, the cone operations related to zero cones and nonnegative cones are simply vector additions and multiplications, 
which have already been parallelized in CUDA. In our implementation we aggregate all zero cones into a single zero cone in our preprocessing step. 
The same holds for nonnegative cones. For the remaining cones, we parallelize within each \emph{family} of cones (\emph{i.e.}, second-order cone, exponential cone, etc.). 
For each family of cones, we allocate a thread to each cone belonging to that family, and execute in parallel. This adheres to the SIMD computing paradigm, as each type of cone has its own set of instructions for updating its scaling matrix. 

We employ our mixed parallel computing strategy in step 4, the scaling matrix update. Recall that the conic constraint $\mcf{K}$ 
can be decomposed as a Cartesian product of $p$ constituent atomic cones $\mcf{K}_1 \times \dots \times \mcf{K}_p$ that are ordered by cone types. 
We assume $\mcf{K}_1$ is a zero cone, $\mcf{K}_2$ is a nonnegative cone, $\mcf{K}_3$ to $\mcf{K}_i$ are second-order cones, 
$\mcf{K}_{i+1}$ to $\mcf{K}_{j}$ are exponential cones, $\mcf{K}_{j+1}$ to $\mcf{K}_l$ are power cones and $\mcf{K}_{l+1}$ to $\mcf{K}_p$ are positive semidefinite cones. Due to our decomposition into constituent atomic cones, 
$H$ is a block-diagonal matrix,
\begin{align*}
	H = \begin{bmatrix}
		H_1 & & \\
		& \ddots & \\
		& & H_p
	\end{bmatrix},
\end{align*}
where $H_1$ and $H_2$ are diagonal matrices and each block $H_t, t \ge 3$ is the scaling matrix corresponding to cone $\mcf{K}_t$ of small dimensionality.

{
	\begin{algorithm}
		\caption{Algorithmic sketch for the scaling matrix $H^k$ update in Algorithm~\ref{alg:IPM}}
		\begin{algorithmic}[1]
			\REQUIRE Current iterate $v^k$, streams $st_{soc}, st_{exp}, st_{pow}, st_{sdp}$.
			\STATE \texttt{function \_update\_H}($v^k$)
			\STATE \qquad Update $H_{\text{zero}}^{k}$  \qquad	// Zero cone ($t = 1$)
			\STATE \qquad Update $H_{\text{nn}}^{k}$	\qquad // Nonnegative cone ($t = 2$)
			\STATE
			\STATE \qquad // Second-order cones
			\STATE \qquad $H_{\text{soc}}^{k} = $ \texttt{\_kernel\_soc\_update\_H<$st_{soc}$>}() // $t = 3$ to $i$
			\STATE
			\STATE \qquad // Exponential cones
			\STATE \qquad $H_{\text{exp}}^{k} = $ \texttt{\_kernel\_exp\_update\_H<$st_{exp}$>}() // $t = i+1$ to $j$
			\STATE
			\STATE \qquad // Power cones
			\STATE \qquad $H_{\text{pow}}^{k} = $ \texttt{\_kernel\_pow\_update\_H<$st_{pow}$>}() // $t = j+1$ to $l$		
			\STATE
			\STATE \qquad // Positive semidefinite cones
			\STATE \qquad $H_{\text{sdp}}^{k} = $ \texttt{\_kernel\_sdp\_update\_H<$st_{sdp}$>}() // $t = l+1$ to $p$	
			\STATE
			\STATE \qquad // Synchronize scaling update
			\STATE \qquad \texttt{CuDeviceSynchronize()}
			\STATE
			\STATE \qquad \texttt{return} $H^{k}$		\qquad // Output the scaling matrix $H^{k}$
			\STATE \texttt{end}
		\end{algorithmic}\label{alg:update-H}
	\end{algorithm}
}
Algorithm~\ref{alg:update-H} illustrates the update of scaling matrix $H^k$ in function \texttt{\_update\_H}($v^k$). 
For the update of $H$ at iteration $k$, we first set diagonal terms of the scaling matrix $H_1^{k}$ for the zero cone to all $0$s, 
and then update the diagonals of the scaling matrix $H_2^{k}$ for the nonnegative cone by element-wise vector division, 
as in the Nesterov-Todd (NT) scaling~\cite{Nesterov98}. 
These computations are already parallelized by CUDA. For conic constraints other than the zero cone and the nonnegative cone, the corresponding parts 
in the scaling matrix $H$ are no longer diagonal, but we observe $H_t$ only requires local information within each constituent cone, and thus we can 
solve for each $H_t$ concurrently, independent to the other cones. 
{We implement kernel functions \texttt{\_kernel\_soc\_update\_H()}, \texttt{\_kernel\_exp\_update\_H()}, 
\texttt{\_kernel\_pow\_update\_H()} and \texttt{\_kernel\_sdp\_update\_H()} for second-order cones, exponential cones, 
power cones and positive semidefinite cones respectively. 
Details regarding the update functions for each type of cone can be found in~\S\ref{appendix:scaling-matrix}.
Each kernel function will process one cone per thread, 
and updating scaling matrices of the same class of cones is executed in parallel. 
Kernel functions for different cone classes are launched independently in separate streams.
Synchronization is required at the end of the scaling matrix update.}

Since the mixed parallel computing strategy is independent of choices of optimization algorithms, it is also applicable for GPU implementations for cone operations 
in first-order operator-splitting conic solvers, like SCS~\cite{scs} and COSMO~\cite{COSMO}.

Note that the SIMD GPU computing structure naturally favors \emph{balanced workloads}, \emph{i.e.}, the workloads in each thread should be similar 
so that the synchronization will not take too much time. The exponential and power cones are 3-dimensional nonsymmetric cones that 
all cone operations have balanced workload among the same class of cones, thus we can parallelize the computation for each cone. 
On the other hand, the second-order cones may vary in dimensionality, which may not mesh well with the SIMD computing paradigm. 
In most use cases for large-scale second-order cone problems, \emph{e.g.}, optimal power flow problems~\cite{PGLIB:2019} and finite-element problems~\cite{Makrodimopoulos2007}, 
the second-order cones are of small (less than 5) dimensionality; thus the effect of this workload imbalance is negligible. Regardless, we support second-order cones of all sizes. 
In the next section, we outline how to process second-order cones of high dimensionality.

\subsection{Dynamic parallelism for second-order cones}\label{subsec:soc-preprocessing}
{
Due to the definition of the barrier function~\eqref{eq:soc-barrier} over a second-order cone \((t, x) \in \mathcal{K}_{q}^n\), 
the computation of residuals involving \(t^2 - \|x\|^2\) requires substantial reduction operations within each cone. 
This introduces an additional level of parallelism, particularly relevant when the dimension of a second-order cone 
reaches thousands or more, as seen in applications like multistage portfolio optimization.
}

{
In such high-dimensional cases, processing second-order cones in a purely thread-wise fashion 
above becomes inefficient, as the residual computation becomes effectively single-threaded. 
A straightforward workaround is to use the dot product operations provided by cuBLAS for computing residuals. 
However, this approach processes each high-dimensional cone sequentially and fails to 
exploit cone-level parallelism 
when multiple second-order cones are present.
To address this, we implement custom reduction operations and leverage \emph{dynamic parallelism} in CUDA, 
which enables additional cone-level parallelism and better utilization of GPU resources.
}

{
Algorithm~\ref{alg:dynamic-parallelism} illustrates how to compute $q$ residuals via dynamic parallelism. 
Suppose we have a vector $x := [x_1; \dots; x_q]$ where $x_i := (t_i, u_i) \in \mathcal{K}_{q}^{n_i}$ 
and residual vector $r \in \Re^{q}$ that stores the residual of cone $i$, \emph{i.e.}, $r[i] = t_i^2 - \|u_i\|^2$. 
We first launch \texttt{\_parent\_kernel\_soc\_residual}($x, r$) for the residual computation. 
We then launch $q$ child kernel functions, \texttt{\_child\_kernel\_soc\_residual} ($x_i, r$), for each thread $i$ ($1 \le i \le q$) 
within the parent kernel function, processing the residual of cone $i$.
The $2$-norm of $u_i$ can be realized by the standard parallel reduction utilizing shared memory efficiently, see Algorithm~\ref{alg:norm2-parallel-reduction}. Finally, we obtain the residual for each cone $i$ and store it back to $r$ in the child function.
}
\begin{algorithm}
	\caption{Dynamic parallelism for the residual computation of $q$ second-order cones}
	\begin{algorithmic}[1]
		\REQUIRE $x = [x_1; x_2; \dots; x_q]$ and $r \in \Re^{q}$, where $x_i = [t_i; u_i] \in \Re^{n_i + 1}$
		\STATE
		\STATE // Parent function
		\STATE \texttt{function \_parent\_kernel\_soc\_residual}($x, r$)
		\STATE \qquad $i \leftarrow$ \texttt{(blockIdx.x -1)* blockDim.x + threadIdx.x}	\quad //Julia is 1-indexed
		\STATE \qquad \textbf{if} ($i \le q$)
		\STATE \qquad \qquad $x_i \leftarrow$ extract cone block $i$ from $x$
		\STATE 
		\STATE \qquad \qquad //compute residual of each cone $i$
		\STATE \qquad \qquad \texttt{\_child\_kernel\_soc\_residual}($x_i, r$)
		\STATE \qquad \textbf{end}
		\STATE \texttt{end}
		\STATE
		\STATE // Child function
		\STATE \texttt{function \_child\_kernel\_soc\_residual}($x_i, r$)
		\STATE \qquad $t_i \leftarrow x_i[1]$ \hfill // Julia: 1-based indexing
		\STATE \qquad $u_i \leftarrow x_i[2:end]$
		\STATE \qquad $\|u\|_2^2 \leftarrow \texttt{\_reduction\_norm2}(u_i, n_i)$
		\STATE \qquad $r[i] \leftarrow t_i^2 - \|u_i\|_2^2$
		\STATE \texttt{end}
	\end{algorithmic}\label{alg:dynamic-parallelism}
\end{algorithm}

{
\subsection{Batched support for positive semidefinite cones}
The positive semidefinite (PSD) cone is a matrix cone that can vary in dimensionality.
Computing scaling matrices and step sizes involves several matrix factorizations, 
such as Cholesky factorization, singular value decomposition and eigenvalue decomposition~\cite{cvxopt}. 
In CuClarabel, we address this by using batched matrix factorizations provided by the cuSOLVER library. 
Hence, we only support a restricted class of SDPs in which all PSD cones have the same dimensionality that is less than or equal to $32$. 
This includes examples from finite element analysis problems~\cite{FEM1,FEM2}, 
where each SDP constraint encodes localized, element-wise or point-wise conditions
—such as stress admissibility or yield criteria—that can be uniformly applied across the mesh.
}

\subsection{Data structures}
\begin{figure}[h]
	\centering
	\begin{tikzpicture}
		% Array of Structures (AoS)
		\node at (2.5, 0) {AoS};
		
		% First structure
		\draw[fill=blue!20] (1, 3) rectangle (4, 4) node[pos=.5, yshift=0.25cm] {$\mathcal{K}_1$};
		\draw[fill=red!20] (1.5, 3) rectangle (2.5, 3.5) node[pos=.5] {$x_1$};
		\draw[fill=green!20] (2.5, 3) rectangle (3.5, 3.5) node[pos=.5] {$y_1$};
		
		% Dots
		\draw[fill=blue!20] (1, 2) rectangle (4, 3) node[pos=.5] {\ldots};
		
		% Last structure
		\draw[fill=blue!20] (1, 1) rectangle (4, 2) node[pos=.5, yshift=0.25cm] {$\mathcal{K}_p$};
		\draw[fill=red!20] (1.5, 1) rectangle (2.5, 1.5) node[pos=.5] {$x_p$};
		\draw[fill=green!20] (2.5, 1) rectangle (3.5, 1.5) node[pos=.5] {$y_p$};
		
		% Structure of Arrays (SoA)
		\node at (8.5, 0) {SoA};
		
		% Whole K structure
		\draw[fill=blue!20] (6.75, 0.5) rectangle (10.25, 4.25) node[pos=.5] {$\mathcal{K}$};
		
		% Array x
		\draw[fill=red!20] (7, 1) rectangle (8, 4);
		\draw[fill=red!20] (7, 2) rectangle (8, 1) node[pos=.5] {$x_p$};
		\draw[fill=red!20] (7, 3) rectangle (8, 2) node[pos=.5] {$\vdots$};
		\draw[fill=red!20] (7, 4) rectangle (8, 3) node[pos=.5] {$x_1$};
		
		% Array y
		\draw[fill=green!20] (9, 1) rectangle (10, 4);
		\draw[fill=green!20] (9, 2) rectangle (10, 1) node[pos=.5] {$y_p$};
		\draw[fill=green!20] (9, 3) rectangle (10, 2) node[pos=.5] {\vdots};
		\draw[fill=green!20] (9, 4) rectangle (10, 3) node[pos=.5] {$y_1$};
		
		% Labels
		\node at (7.5, 0.75) {$x$};
		\node at (9.5, 0.75) {$y$};
		
		% Bounding box
		\draw[draw=none] (-1, -1) rectangle (12, 6); % Invisible bounding box for better spacing
	\end{tikzpicture}
	\caption{Illustration of AoS (CPU) and SoA (GPU) data structure}
	\label{fig:AoS-vs-SoA}
\end{figure}
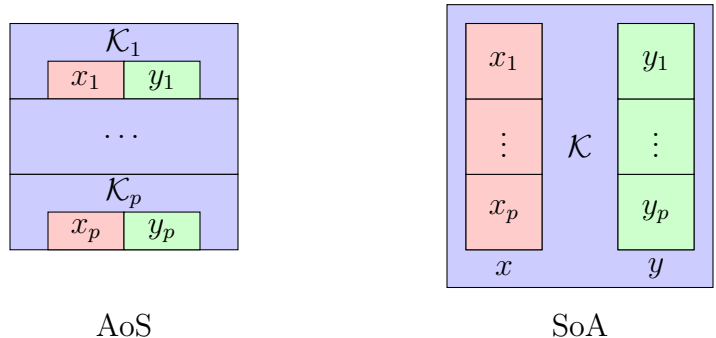

Since data parallelism is required for parallel execution on GPUs, we manage data for each cone as a structure of arrays (SoA) in our GPU implementation, 
in contrast to the existing CPU counterpart which uses an array of structures (AoS). In other words, instead of creating structs for $\mcf K_i, i=1,\dots,p$ 
and storing pointers to each struct in an array, we concatenate the same type of local variable from different cones into a global variable and 
then store it in a global struct for $\mcf K$, which is illustrated in Figure~\ref{fig:AoS-vs-SoA}.

In addition, indexing cones in order in the setup phase can simplify the implementation of mixed parallel computing strategy. 
We reorder the input cones such that memory is coalescing for the same class of cones, which can accelerate computation on GPUs.
Matrices are stored in the compressed sparse row (CSR) format. Instead of storing only the triangular part of a square matrix as in the CPU-based Clarabel, 
we store the full matrix for more efficient multiplication on the GPU.

\subsection{Solving linear systems}
Most of the computation time for our interior point method is spent in factorizing the matrix $K$ and in the three backsolve operations in~\eqref{eqn:linsys_2x2}. 
Our GPU implementation also leverages the power of the newly released sparse linear system solver cuDSS~\cite{CUDSS} for the  $LDL^T$ factorization and backsolve operations. 
{We implement a self-contained iterative refinement originated from~\cite{Nocedal06} to increase the numerical stability of the backsolve operation.}

Currently, a GPU has more computing cores for \textsc{Float32} than \textsc{Float64} and hence better performance for parallel algorithms. 
Also, \textsc{Float32} requires less memory and takes less time for the same computation than \textsc{Float64}. 
However, lowering the precision will introduce numerical instability for solving linear systems. Thus, we employ a \emph{mixed precision} solve in our matrix solves. 
Specifically, we use mixed precision for data within the iterative refinement step:

\begin{enumerate}
	\item Solve for the residual at step $i$: $r_i = b - Kx_{i-1}$, where $x_{i-1}$ is our guess for iteration $i-1$. (Full precision)
	\item Solve the linear system under a regularization parameter $(K+ \delta I)\Delta_i = r_i$. (Lower precision)
	\item Take $x_i = x_{i-1} + \Delta_i$. (Full precision)
\end{enumerate} 

In all our matrix solves, given matrix $K$ and right-hand side $b$, we factorize the lower precision copy of $K$, 
and then solve the linear system in this lower precision. 
{
To enhance numerical stability during matrix factorization, a regularization term $\delta I$ with 
$$
\delta := \delta_s + \delta_d \max_i |D_{ii}|,
$$
is added to $K$, where $\delta_s$ is a static regularization and $\delta_d \max_i |D_{ii}|$ provides dynamic scaling based on the magnitude of the diagonal entries of $D$.
}
To offset the regularization effect and rounding errors from the lower precision in backsolves, 
we solve in full precision for the other steps, \emph{i.e.}, we compute the residual $r_i$ and save the update $\{x_i\}$ in full precision. {When we change the factorization data type from \textsc{Float32} to \textsc{Float64}, it reduces down to the standard iterative refinement as employed in standard conic optimization solvers~\cite{Clarabel, ecos}.}

{We apply iterative refinement until either (1) we reach a pre-specified number 
of max steps $T_{\max}$ for iterative refinement, or (2) the $\ell_\infty$ norm 
of $b - Kx_i$ reaches a certain threshold, 
specifically $\norm{b - Kx_i}_\infty \leq t_\text{abs} + t_\text{rel}  \norm{b}_\infty$. 
We set both $t_\text{abs}$ (absolute tolerance) and $t_\text{rel}$ (relative tolerance) 
to $10^{-12}$, and the maximum number of steps $T_{\max} = 10$. 
We set $\delta_s$ to the square root of machine precision, 
\emph{i.e.}, $\sqrt{\epsilon}\approx 3.45e^{-4}$, for mixed precision and $\delta_s=1e^{-8}$ for full precision, and we set $\delta_d$ to the square of machine precision, \emph{i.e.}, $\epsilon^2$.}

Although the mixed precision for the iterative refinement can improve the numerical stability for solving linear system in lower precision, 
it is to be expected that the mixed precision may take longer time to converge or fail in cases where the matrix $K$ is extremely ill-conditioned. 
Caution is required when using mixed precision for numerically hard problems, \emph{e.g.}, conic programs with exponential cones. {Hence, we provide the mixed precision as an optional choice and use \textsc{Float64} by default.}

\section{Numerical experiments}\label{sec:numerical}
We have benchmarked our Julia GPU implementation of the Clarabel solver, {denoted as \emph{ClarabelGPU},} against the state-of-the-art commercial interior-point solvers MOSEK~\cite{mosek} and Gurobi~\cite{gurobi}. 
We also include our Rust CPU implementation of Clarabel solver with the 3rd-party multithreaded supernodal LDL factorization method
in the \emph{faer-rs} package~\cite{faer}, {denoted as \emph{ClarabelRs}. Note that MOSEK and Gurobi already utilize multithreaded linear system solvers}. We include benchmark results for several classes of problems including quadratic programming (QP), 
second-order cone programming (SOCP) and exponential cone programming. All benchmarks are performed using the default settings for each solver, 
with pre-solve disabled where applicable to ensure equivalent problem-solving conditions. 
No additional iteration limits are imposed beyond each solver's internal defaults.	
{
We set $\epsilon_{\text{feas}}=1e^{-6}$ for all solvers in termination check.
}
All experiments were carried out on a workstation with Intel(R) Xeon(R) w9-3475X CPU @ 4.8 GHz 
with 256 GB RAM and NVIDIA GeForce RTX 4090 24GB GPU. 
All benchmarks tests are scripted in Julia and access solver interfaces via JuMP [82]. 
We use Rust compiler version 1.76.0 and Julia version 1.10.2. {Before recording benchmark tests, we run a few example solves; 
this offsets the just-in-time (JIT) precompilation overhead in Julia.}

\subsection{Benchmarking metrics}
We choose the same benchmarking tests as used in the Clarabel solver~\cite{Clarabel}, and compare our results using metrics that are commonly used when comparing computational time 
across different solvers. For a set of $N$ test problems, we define the \emph{shifted geometric mean} $g_s$ as 
\[
g_s \eqdef \left[\prod_{p=1}^N  (t_{p,s} + k) \right]^{\frac{1}{N}} - k,
\]
where $t_{p,s}$ is the time in seconds for solver $s$ to solve problem $p$, and $k = 1$ is the shift. The normalized shifted geometric mean is then defined as 
\[
r_s \eqdef \frac{g_s}{\min_{s'} g_{s'}}.
\]
Note that the solver with the lowest shifted geometric mean solve time has a normalized score of~1. We assign a solve time $t_{p,s}$ equal to the maximum allowable solve time 
if solver $s$ fails to solve the problem $p$.

The \emph{relative performance ratio} for a solver $s$ and a problem $p$ is defined as
\[
u_{p,s} = \frac{t_{p,s}}{\min_{s'} t_{p,s'}}.
\]
The \emph{relative performance profile} denotes the fraction of problems solved by solver $s$ within a factor $\tau$ of the solve time of the best solver, 
which is defined as $f^r_s : \Re_+ \mapsto [0,1]$  
\[
f^r_s(\tau) \eqdef \frac{1}{N} \sum_p \mathcal{I}_{\le \tau}(u_{p,s}),
\]
where $\mathcal{I}_{\le \tau}(u) = 1$ if $\tau \le u$ and $\mathcal{I}_{\le \tau}(u) = 0$ otherwise. {For $\tau = 1$, $f^r_s(1)$ denotes the ratio of problems where the solver $s$ performs the best, and $\sum_{s} f^r_s(1)$ should be equal to 1 if every problem is solvable for at least one solver.} We also compute the \emph{absolute performance profile} 
$f^a_s : \Re_+ \mapsto [0,1]$, which denotes the fraction of problems solved by solver $s$ within $\tau$ seconds and is defined as
\[
f^a_s(\tau) \eqdef \frac{1}{N} \sum_p \mathcal{I}_{\le \tau}(t_{p,s}).
\]

{Our numerical experiments highlight the performance gains achieved by the GPU 
implementation of Clarabel on a diverse set of conic optimization problems. We record \emph{total time} by default, the sum of the \emph{setup time}, including data equilibration and solver initialization, and the \emph{solve time}, which is the running time for the algorithm underlying a solver.} 

\subsection{Quadratic programming}\label{subsec:qp}
We first present benchmark results for QPs. Note that in this setting, the set $\mathcal{K}$ in~\eqref{eqn:primal} is restricted to the composition of zero cones, 
\emph{i.e.}, linear equality constraints, and nonnegative cones, \emph{i.e.}, linear inequality constraints. We consider two classes of problems, the portfolio optimization problem 
and the Huber fitting problem.

Portfolio optimization, a problem arising in quantitative finance, aims to allocate assets in a manner that maximizes expected return while keeping risk under control. 
We can formulate it as 
$$
\begin{array}{ll}
\text{maximize}    & \mu^T x-\gamma x^T \Sigma x \\
\text {subject to} & \mathbf{1}^T x=1, \\
                   & x \geq 0,
\end{array}
$$
where $x \in \Re^n$ (the variable) represents the ratio of allocated assets, $\mu \in \Re^n$ is the vector of expected returns, $\gamma>0$ is the risk-aversion parameter, 
and $\Sigma \in \mathbf{S}_{+}^n$ the risk covariance matrix which is of the form $\Sigma \eqdef  FF^T + D$ with $F \in \Re^{p \times n}$  and $D \in \Re^{p \times p}$ diagonal. 
We set the rank $p$ to the integer closest to $0.1n$, and vary $n$ from $5000$ to $25000$.

Huber fitting is a version of robust least squares. For a given matrix $A \in \Re^{m \times n}$ and vector $b \in \Re^m$, we replace the least squares loss function with the Huber loss. 
The Huber loss makes the penalty incurred by larger points linear instead of quadratic, thus outliers have a smaller effect on the resulting estimator. Precisely, the problem is stated as:
$$
\text{minimize} \ \sum_{i=1}^m \phi_{h}\left(a_i^T x-b_i\right)
$$
where $a_i^T$ is the $i$-th row of $A$ and the Huber loss $\phi_{h}: \Re \rightarrow \Re$ is defined as
$$
\phi_{h}(t)\eqdef \begin{cases}t^2 & |t| \leq T \\ T(2|t|-T) & \text { otherwise }\end{cases}.
$$
We set $m$ to the nearest integer of $1.5n$ and vary value of $n$ from $5000$ to $25000$.
{
	
	\captionsetup{labelfont=bf}
	\centering
	\renewcommand{\detailtablecaption}{\bf Benchmarking for quadratic programming}

	\tiny
\begin{longtable}{||l||cccc||cccc||}
	\caption{\detailtablecaption}
	\label{table:qp}
	\\
	& \multicolumn{4}{c||}{\underline{iterations}}& \multicolumn{4}{c||}{\underline{total time (s)}}\\[2ex] 
	Problem & ClarabelGPU & Mosek* & ClarabelRs & Gurobi & ClarabelGPU & Mosek* & ClarabelRs & Gurobi\\[1ex]
	\hline
	\endhead
	\sc{portfolio\_optimization\_5000}  &   15 &  \winner  9 &   15 &   16 &  \winner  0.372 &   0.936 &   7.36 &   0.513\\ 
	\sc{portfolio\_optimization\_10000} &   18 &  \winner  10 &   18 &   17 &  \winner  1.42 &   4.24 &   48.8 &    2.8\\ 
	\sc{portfolio\_optimization\_15000}  &   18 &  \winner  11 &   18 &   17 &  \winner  3.55 &   12.7 &    147 &   7.52\\ 
	\sc{portfolio\_optimization\_20000}  &   17 &  \winner  12 &   17 &   17 &  \winner  6.66 &   28.5 &    307 &   14.7\\ 
	\sc{portfolio\_optimization\_25000}  &   18 &  \winner  11 &   18 &   18 &  \winner  11.8 &   49.3 &    549 &   27.2\\ 
	\sc{huber\_fitting\_5000}  &  \winner  8 &   11 &  \winner  8 &   9 &  \winner  4.81 &   13.5 &   65.6 &   48.2\\ 
	\sc{huber\_fitting\_10000}  &  \winner  8 &   11 &  \winner  8 &   9 &  \winner  37.1 &   55.1 &    612 &    276\\ 
	\sc{huber\_fitting\_15000}  &  \winner  8 &   13 &  \winner  8 &   10 &  \winner   122 &    289 &   2.1e+03 &    747\\ 
	\sc{huber\_fitting\_20000} &  \winner  8 &   13 &   - &   10 &  \winner   283 &    807 &   - &   1.76e+03\\ 
	\sc{huber\_fitting\_25000} &  \winner  8 &   15 &   - &   10 &  \winner   547 &   1.83e+03 &   - &   1.24e+03\\ 
\end{longtable}
	
	%\label{table:#1:detail}

}
Results for large QP tests are shown in Table~\ref{table:qp}. We benchmark ten different examples from two classes above and set the time limit to $1$h. 
We compare our GPU implementation ClarabelGPU with ClarabelRs 
and two commercial solvers, Gurobi and MOSEK. 
ClarabelGPU is the fastest solver on these problems, 
and it has the lowest per-iteration time for almost all examples. 
Since most of time of an interior point solver is spent on factorizing and solving a linear system in QPs, we can say that ClarabelGPU benefits from the use of the cuDSS linear system solver 
and it is more than 2 times faster than Gurobi, about 4 times faster than MOSEK and 10x times faster than the existing Rust implementation with the multithreaded \emph{faer-rs} linear system solver. 

\subsection{Second-order cone programming}
We next consider the second-order cone relaxations of optimal power flow problems~\cite{PowerModels:SOCcase} from the IEEE PLS PGLib-OPF benchmark library~\cite{PGLIB:2019}, 
using the \emph{PowerModels.jl} package~\cite{PowerModels:2018} for modeling convenience. Note that only second-order cones of dimensionality $3$ or $4$ are used in these second-order cone relaxations, 
which satisfies our assumption in \S\ref{subsec:kernel-operations}: that the dimensionality of each cone is very small.

We compare our GPU implementation with the CPU-based Clarabel solver and MOSEK solver. We also include results of the MOSEK solver with pre-solve for comparison, 
which is denoted as \texttt{Mosek*} in the plots. The maximum termination time is again set to $1$h. We benchmark the second-order cone relaxations of 120 problems from the PGLib-OPF library, 
where the number of second-order cones exceed $2000$.

\begin{figure}[H]
		\centering
		\caption{\bf Performance profiles for the large OPF SOCPs problem set}
		\label{fig:gpu_opf_large_socp}
		\begin{subfigure}[b]{0.49\textwidth}
			\centering
			{\includegraphics[width=\textwidth]{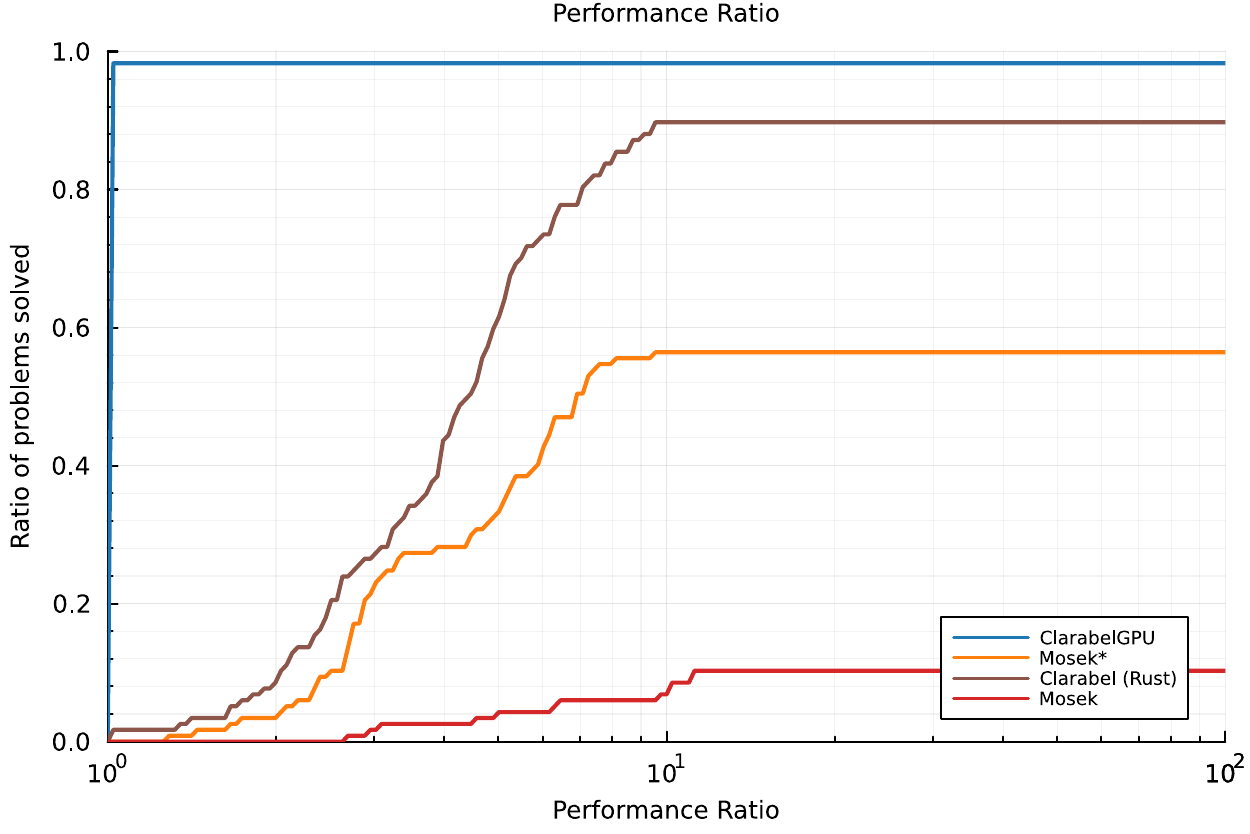}}
			\caption{Relative performance profile}
			\label{fig:gpu_opf_large_socp:relative}
		\end{subfigure}
		\hfill
		\begin{subfigure}[b]{0.49\textwidth}
			\centering
			{\includegraphics[width=\textwidth]{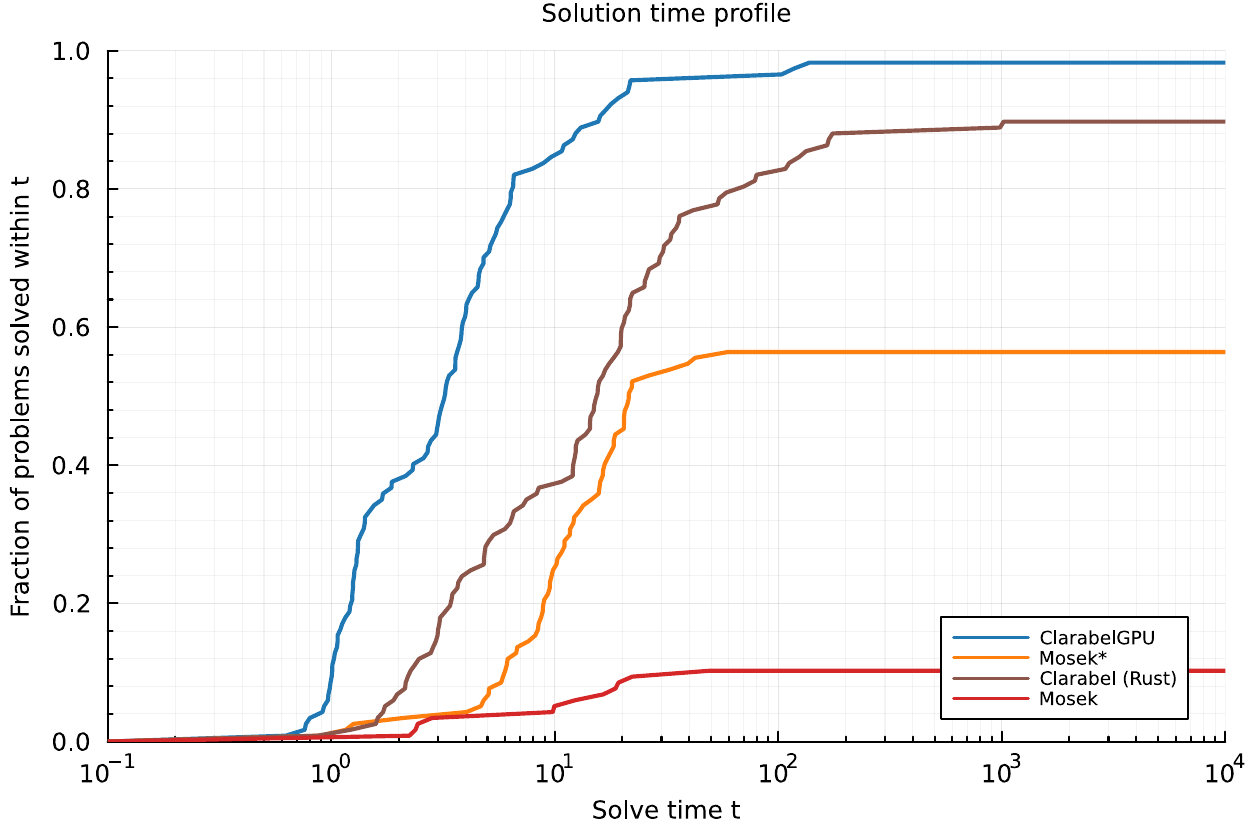}}
			\caption{Absolute performance profile}
			\label{fig:gpu_opf_large_socp:absolute}
		\end{subfigure}
		\begin{subfigure}{1\textwidth}
			\centering
			\footnotesize
			\begin{tabular}{llcccc}
  \hline
   &  & \textbf{ClarabelGPU} & \textbf{Mosek*} & \textbf{ClarabelRs} & \textbf{Mosek} \\\hline
  Shifted GM & Full Acc. & 1.0 & 12.12 & 4.59 & 53.65 \\
   & Low Acc. & 1.0 & 4.16 & 3.91 & 7.29 \\
  Failure Rate (\%) & Full Acc. & 1.7 & 43.6 & 10.3 & 89.7 \\
   & Low Acc. & 0.0 & 5.1 & 1.7 & 0.9 \\\hline
\end{tabular}

			\caption{Benchmark timings as shifted geometric mean and failure rates}
		\end{subfigure}
	\end{figure}

Results for these problems are shown in Figure~\ref{fig:gpu_opf_large_socp}. Both ClarabelGPU and ClarabelRs are faster and more numerically stable than MOSEK even with the presolve step.
Moreover, the GPU implementation can solve 118 out of 120 examples within an hour, while the Rust version of Clarabel can solve 104 out of 120 examples within the same time limit. 
In contrast, MOSEK with presolve fails on about 40\% of the optimal flow problems, and the one without presolve fails on over 90\% of the problems. 

Overall, the GPU solver is several times faster than Rust-based CPU solver. {This speedup causes ClarabelGPU to achieve a lower failure rate than ClarabelRs as it could solve more problems in under $10^4$ seconds each}. However, note that  ClarabelGPU fails on two examples that ClarabelRs successfully solves.
This highlights the different numerical performance of the linear system solvers between cuDSS and \emph{faer-rs}.

	\captionsetup{labelfont=bf}

{We also consider the following multistage portfolio optimization problem:}
\begin{equation}\label{eq:multistage_socp}
\begin{aligned}
	\text{minimize} \quad 
	& \sum_{t=1}^{T} \left( -\mu_t^\top x_t + c \cdot \mathbf{1}^\top z_t + \gamma \cdot r_t \right) \\
	\text{subject to} \quad
%	& x_t \in [0, 0.1]^n, \quad
%	  y_t \in [0, 0.1]^k, \quad
%	  z_t \in \mathbb{R}^n, \quad
%	  r_t \ge 0, \quad \forall t, \\
	& z_1 \ge x_1 - x_0, \quad z_1 \ge x_0 - x_1, \\
	& z_t \ge x_t - x_{t-1}, \quad z_t \ge x_{t-1} - x_t, \quad \forall t = 2, \ldots, T, \\
	& \mathbf{1}^\top x_1 = d + \mathbf{1}^\top x_0, \quad
	  \mathbf{1}^\top x_t = \mathbf{1}^\top x_{t-1}, \quad \forall t = 2, \ldots, T, \\
	& y_t = F_t x_t, \quad \forall t = 1, \ldots, T, \\
	& \left[
	\begin{array}{c}
		r_t \\
		U y_t \\
		D_{\text{sqrt}} \odot x_t
	\end{array}
	\right] \in \mathcal{K}_{q}^{n+k+1}, \quad \forall t = 1, \ldots, T,
\end{aligned}
\end{equation}
{with respect to $x$, $y$, $z$, and $r$. This problem is formulated over $T$ time periods, allocating wealth across $n$ assets while managing exposure to $k$ underlying factors. Here, $x_t$ represents asset allocations, $y_t$ represents factor exposures, $z_t$ represents trade volumes, and $r_t$ is a risk proxy. We impose additional constraints $x_t \in [0, 0.1]^n$, $y_t \in [0, 0.1]^k$, $z_t \in \mathbb{R}^n$, and $r_t \geq 0$. Our objective is a minimization over the trade-off between negative expected returns, linear transaction costs, and a second-order cone-based risk penalty scaled by risk aversion $\gamma$. The transaction constraints ensure that trade volumes \( {z}_t \) correctly reflect changes in portfolio weights. Budget constraints enforce capital conservation: initial capital changes by a known inflow \( d \), and total capital remains constant thereafter. Factor exposure is modeled via linear mappings \( {y}_t = F_t {x}_t \). Risk $r_t$ is quantified through a second-order cone (SOC) constraint that captures both systematic risk \( U y_t\) and idiosyncratic risk \( D_{\text{sqrt}}\odot x_t \) across each time period.}

{
We benchmark multistage portfolio optimization problems with $n=5000, k=50$ and varying horizon $T$ in Table~\ref{table:parametric_programming_large_socp}; 
ClarabelGPU achieved up to a $3$x speedup over MOSEK and $40$x over ClarabelRs.}

{
	
	\captionsetup{labelfont=bf}
	\centering
	\renewcommand{\detailtablecaption}{\bf Multi-stage portfolio optimization with varying horizons}

	\scriptsize
\begin{longtable}{||l||cccc||cccc||}
	\caption{\detailtablecaption}
	\label{table:parametric_programming_large_socp}
	\\
	& \multicolumn{4}{c||}{\underline{iterations}}& \multicolumn{4}{c||}{\underline{total time (s)}}\\[2ex] 
	Horizon  & ClarabelGPU & Mosek* & ClarabelRs & Mosek & ClarabelGPU & Mosek* & ClarabelRs & Mosek\\[1ex]
	\hline
	\endhead
	\sc{5} &   21 &   15 &   21 &  \winner  14 &  \winner  0.882 &   1.24 &   12.1 &   1.65\\ 
	\sc{10} &   22 &  \winner  15 &   22 &  \winner  15 &  \winner  1.79 &   3.87 &   53.7 &   4.63\\ 
	\sc{15} &   22 &   14 &   22 &  \winner  12 &  \winner   2.8 &   7.49 &   99.5 &   7.42\\ 
	\sc{20} &   22 &   15 &   22 &  \winner  11 &  \winner  3.91 &   12.5 &    167 &   13.2\\ 
	\sc{25} &   24 &   18 &   24 &  \winner  14 &  \winner  9.82 &   20.2 &    245 &   19.8\\ 
	\sc{30} &   23 &  \winner  12 &   23 &   13 &  \winner  11.3 &   22.1 &    312 &     28\\ 
\end{longtable}
	
	%\label{table:#1:detail}

}
\subsection{Exponential cone programming}
{For exponential programming, we benchmark the entropy maximization problems~\cite{scs} with varying dimensionality. } The entropy maximization problem aims to maximize entropy over a probability distribution given a set of $m$ linear inequality constraints, 
which can be interpreted as bounds on the expectations of arbitrary functions. The problem is formulated as 
$$
\begin{array}{ll}
\text{maximize} & -\sum_{i=1}^n x_i \log x_i \\
\text {subject to} & 1^T x=1, \\
                   & A x \leq b,
\end{array}
$$
Each element of $A$ is generated from the distribution $A_{i j} \sim \mathcal{N}(0, n)$. Then, we set $b=A v /1^T v$ 
where $v \in \Re^n$ is generated randomly from $v_i \sim U[0,1]$, which ensures the problem is always feasible. 
We set $m$ to the nearest integer of $0.5n$ and vary value of $n$ from $2000$ to $10000$.

{

	\captionsetup{labelfont=bf}
	\centering
	\renewcommand{\detailtablecaption}{\bf Benchmarking for exponential cone programming}

	\scriptsize
\begin{longtable}{||l||cccc||cccc||}
	\caption{\detailtablecaption}
	\label{table:exponential_cone}
	\\
	& \multicolumn{4}{c||}{\underline{iterations}}& \multicolumn{4}{c||}{\underline{total time (s)}}\\[2ex] 
	Problem & ClarabelGPU & Mosek* & ClarabelRs & Mosek & ClarabelGPU & Mosek* & ClarabelRs & Mosek\\[1ex]
	\hline
	\endhead
	\sc{entropy\_max\_2000} &   20 &   15 &   20 &  \winner  10 &  \winner  0.476 &   0.794 &   2.84 &   0.609\\ 
	\sc{entropy\_max\_4000} &   20 &   16 &   20 &  \winner  9 &  \winner   1.4 &   3.02 &     14 &   2.28\\ 
	\sc{entropy\_max\_6000} &   21 &   16 &   21 &  \winner  9 &  \winner  3.16 &   8.85 &   90.5 &   5.89\\ 
	\sc{entropy\_max\_8000} &   21 &   16 &   21 &  \winner  10 &  \winner  5.98 &   17.7 &    104 &   12.1\\ 
	\sc{entropy\_max\_10000} &   21 &   16 &   23 &  \winner  10 &  \winner  19.2 &   31.6 &    239 &   36.5
\end{longtable}
	
	%\label{table:#1:detail}

}

The benchmark results for these two problems with varying dimensionality are shown in Table~\ref{table:exponential_cone}. 
Though the GPU acceleration of Clarabel is offset by nearly doubled number of iterations compared to MOSEK, 
we can still achieve more than 2 times of acceleration for the overall time. That is to say we can possibly achieve more acceleration in ClarabelGPU 
if we can improve the numerical stability to the same level of MOSEK on exponential cone programs. 
We find ClarabelGPU can benefit from GPU computation up to 10x times faster compared to ClarabelRs.

\subsection{Semidefinite programming}
{
Our results in Table~\ref{table:FEM_SDP} show that ClarabelGPU achieves a speedup of $1.5\times$ to $4\times$ over MOSEK, both with and without presolve, and up to $10\times$ speedup compared to ClarabelRs. ClarabelGPU can only solve FELA\_SDP\_9263 to $1e^{-5}$ precision and is less numerical stable than Mosek* under the default setting, but increasing the static regularization $\delta_s$ to $1e^{-7}$ can fix the numerical issue.
}

{

	\captionsetup{labelfont=bf}
	\centering
	\renewcommand{\detailtablecaption}{\bf FEM with $3\times3$-dimensional SDP constraints}

	\scriptsize
\begin{longtable}{||l||cccc||cccc||}
	\caption{\detailtablecaption}
	\label{table:FEM_SDP}
	\\
	& \multicolumn{4}{c||}{\underline{iterations}}& \multicolumn{4}{c||}{\underline{total time (s)}}\\[2ex] 
	Problem  & ClarabelGPU & Mosek* & ClarabelRs & Mosek & ClarabelGPU & Mosek* & ClarabelRs & Mosek\\[1ex]
	\hline
	\endhead
	\sc{FELA\_SDP\_1048} &   23 &   24 &   - &  \winner  22 &  \winner  2.41 &   6.75 &   - &   2.84\\ 
	\sc{OBEFM\_SDP\_1048} &   19 &  \winner  13 &   19 &   15 &  \winner  2.48 &   5.98 &   30.6 &   4.04\\ 
	\sc{FELA\_SDP\_4444} &  \winner  30 &   31 &   - &   38 &  \winner    11 &   40.8 &   - &   19.4\\ 
	\sc{OBEFM\_SDP\_4444} &   20 &  \winner  13 &   20 &   - &  \winner  11.5 &   35.3 &    144 &   -\\ 
	\sc{FELA\_SDP\_9263} &   - &   33 &  \winner  31 &   34 &   - &   55.1 &    167 &  \winner  41.5\\ 
	\sc{OBEFM\_SDP\_9263} &   21 &  \winner  14 &   21 &   - &  \winner  25.4 &   38.1 &    323 &   -\\ 
\end{longtable}
	
	%\label{table:#1:detail}

}

\subsection{Parametric programming}
{
Clarabel supports updating the coefficients $P,A,q,b$ without reinitializing the solver object. 
This feature is particularly advantageous for parametric programming, where the solver setup is required only once, 
and the symbolic factorization structure can be reused efficiently in subsequent solves~\cite{MPT3}. 
Such capability is especially useful in applications like model predictive control and portfolio optimization.
}

{
For the multistage portfolio optimization problem in~\eqref{eq:multistage_socp}, 
we also split the total time into two parts, the setup time and the solve time, 
and report the ratio of them respectively in Table~\ref{table:time_ratio}. Notably, 
the setup time for ClarabelGPU constitutes a larger proportion of the total time compared to its CPU counterpart. 
This suggests that, when the setup time dominates the total computational time, 
higher acceleration ratios can be achieved on a GPU 
when we need to solve multistage portfolio optimization multiple times with different parameters.
}

{{
	\captionsetup{labelfont=bf}
	\centering
	\renewcommand{\detailtablecaption}{\bf Time ratio between setup and solve time on multistage portfolio optimization}

	\scriptsize
\begin{longtable}{||l||cc||cc||}
	\caption{\detailtablecaption}
	\label{table:time_ratio}
	\\
	& \multicolumn{2}{c||}{\underline{ClarabelGPU}}& \multicolumn{2}{c||}{\underline{Clarabel}}\\[2ex] 
	Horizon & setup time & solve time & setup time & solve time\\[1ex]
	\hline
	\endhead
	\sc{5} & 0.648 & 0.352 & 0.418 & 0.582\\ 
	\sc{10} & 0.654 & 0.346 & 0.36 & 0.64\\ 
	\sc{15} & 0.662 & 0.338 & 0.321 & 0.679\\  
	\sc{20} & 0.651 & 0.349 & 0.284 & 0.716\\ 
	\sc{25} & 0.475 & 0.525 & 0.242 & 0.758\\
	\sc{30} & 0.566 & 0.434 & 0.223 & 0.777\\ 
\end{longtable}
	
	%\label{table:#1:detail}
}}

\subsection{Mixed precision}
In the implementation of the mixed precision setting, as described we only set the data type of the cuDSS linear system solver to \textsc{Float32}. 
Since the use of mixed precision accelerates the numerical factorization rather than the symbolic factorization within a factorization method, 
we only record the computational time for an interior point method without the setup time. 

We test the mixed precision setting on QPs, including portfolio optimization and Huber fitting problems from \S\ref{subsec:qp}, for high accuracy level $\epsilon_{\text{feas}}=1e^{-8}$. 
We compare it with the standard GPU implementation of \textsc{Float64} data type. 
{The results in Table~\ref{table:mixed_large_qp} demonstrate that employing a mixed-precision strategy 
can reduce solve time by up to a factor of 2 when the matrix $A$ becomes dense. 
However, we note that mixed precision is less numerically stable when solving an ill-conditioned KKT system~\eqref{eqn:linsys_2x2}. This can lead to an increased number of iterative refinement steps and longer solve times compared to full precision, or even failure to converge to the desired tolerance.
}

{
	\captionsetup{labelfont=bf}
	\centering
	\renewcommand{\detailtablecaption}{\bf Solve times and iteration counts for the mixed precision QP test}

	\scriptsize
\begin{longtable}{||l||cc||cc||}
\caption{\detailtablecaption}
\label{table:mixed_large_qp}
\\
& \multicolumn{2}{c||}{\underline{iterations}} & \multicolumn{2}{c||}{\underline{solve time (s)}}\\[2ex] 
Problem & Full & Mixed & Full & Mixed \\[1ex]
\hline
\endhead
\sc{portfolio\_optimization\_n\_5000} & 17 & 18 &  \winner 0.18 & 0.29\\ 
\sc{portfolio\_optimization\_n\_10000} & 19 & 19 &  0.70 & \winner 0.57\\ 
\sc{portfolio\_optimization\_n\_15000} & 19 & 19 &  1.8 & \winner 1.31\\ 
\sc{portfolio\_optimization\_n\_20000} & 19 & 23 &  3.74 & \winner 3.26\\ 
\sc{portfolio\_optimization\_n\_25000} & 19 & 21 &  6.79 & \winner 5.58\\ 
\sc{huber\_fitting\_n\_5000} & 9 & 9 &  3.6 & \winner 1.85\\ 
\sc{huber\_fitting\_n\_10000} & 10 & 10 &  36.9 & \winner 19.9\\ 
\sc{huber\_fitting\_n\_15000} & 10 & 10&  128 & \winner 70.5\\ 
\sc{huber\_fitting\_n\_20000} & 9 & 9 &  279 &  \winner 156\\
\sc{huber\_fitting\_n\_25000} & 10 & 10 &   598 &  \winner 334\\  
\end{longtable}

	%\label{table:#1:detail}
}

\section{Conclusion}
{
We have developed a GPU interior point solver for conic optimization.\footnote{\url{https://github.com/oxfordcontrol/Clarabel.jl/tree/CuClarabel}}
In our implementation, we propose a mixed parallel computing strategy to process linear constraints with second-order cone, exponential cone, power cone and semidefinite cone constraints. 
Our GPU solver shows several times acceleration compared to state-of-the-art CPU conic 
solvers on many problems to high precision, such as QPs, SOCPs, exponential cone programs and SDPs.
}

{
Future research directions include extending support to general SDPs with PSD cones of varying dimensionalities. 
The proposed mixed parallel computing strategy for GPU implementation is also applicable to conic solvers 
based on first-order operator-splitting methods. 
This approach could improve GPU utilization, especially when sufficient computational resources 
are available to handle different cone classes in parallel. Additional performance gains
may be achieved through kernel fusion, the use of CUDA graphs to reduce kernel launch overhead and overlapping more independent computation within interior-point methods.
}

\section*{Acknowledgments}
{Yuwen Chen was supported by the EPSRC Impact Acceleration Account Fund and the National Centre of Competence in Research (NCCR) Automation, funded by the Swiss National Science Foundation (SNSF)}.
Parth Nobel was supported in part by the National Science Foundation
Graduate Research Fellowship Program under Grant No. DGE-1656518. Any
opinions, findings, and conclusions or recommendations expressed in
this material are those of the author(s) and do not necessarily
reflect the views of the National Science Foundation.
The work of Stephen Boyd was supported by ACCESS (AI Chip Center for Emerging Smart Systems), sponsored by InnoHK funding, Hong Kong SAR.

\bibliographystyle{unsrt}
\bibliography{reference}

\appendix

\clearpage

\section{Scaling matrices}\label{appendix:scaling-matrix}
The scaling matrix $H$ in~\eqref{eqn:linsys_4x4} is the Jacobian from the linearization of the central path~\eqref{eqn:central-path}. We set $H = 0$ for the zero cone and choose the NT scaling~\cite{Nesterov98} for nonnegative and second-order cones, which are both symmetric cones. The NT scaling method exploits the self-scaled property of symmetric cone $\mathcal{K}$ where exists a unique scaling point $w \in \mathcal{K}$ satisfying
\begin{equation*}
	H(w)s = z.
\end{equation*}
The matrix $H(w)$ can be factorized as $H^{-1}(w) = W^T W$, and we set $H = H^{-1}(w)$ in~\eqref{eqn:linsys_4x4}. The factors $w,W$ are then computed following~\cite{cvxopt}.

For exponential and power cones that are not symmetric, the central path is defined by the set of point satisfying
\begin{align*}
	Hz = s, \quad H \nabla f^*(s) = \nabla f(z),
\end{align*}
where $f^*$ is the conjugate function of $f$, and the symmetric scaling from~\cite{Tuncel01,Dahl21} is implemented. We define \emph{shadow iterates} as
\begin{align*}
	\tilde{z} \eqdef -\nabla f(s), \quad \tilde{s}\eqdef - \nabla f^*(z),
\end{align*}
with $\tilde{\mu} = \langle \tilde{s}, \tilde{z} \rangle/\nu$.
A scaling matrix $ H $ is chosen to be the rank-4 Broyden-Fletcher-Goldfarb-Shanno (BFGS) update as in a quasi-Newton method,
\begin{align*}
	H \eqdef H_{\mathrm{BFGS}} \eqdef Z(Z^T S)^{-1}Z^T + H_a - H_aS(S^T H_aS)^{-1}S^T H_a,
\end{align*}
where $Z \eqdef [z, \tilde{z}], S \eqdef [s, \tilde{s}]$ and $H_a = \mu \nabla^2 f(z)$ is an estimate of the Hessian as in~\cite{Dahl21}.

\section{Parallel reduction for norm-$2$ computation}\label{appendix:dynamic-parallelism}
\begin{algorithm}
	\caption{Parallel Reduction to Compute $\|u\|^2$}
	\begin{algorithmic}[1]
		\REQUIRE $u \in \Re^n$, accessible to all threads in a block
		\ENSURE Scalar $\texttt{norm2} = \sum_{j=1}^{n} u_j^2$
		\STATE
		\STATE \texttt{inline function \_reduction\_norm2}(u, n)
		\STATE \qquad Allocate \texttt{shared\_mem[1\,..\,T]} = 0 \hfill // T = threads per block
		\STATE \qquad $j \leftarrow$ \texttt{threadIdx.x}
		\STATE \qquad \textbf{while} $j \le n$ \textbf{then}
		\STATE \qquad \qquad \texttt{shared\_mem[j]} $\leftarrow \texttt{shared\_mem[j]} + u[j]^2$ \hfill // Julia is 1-indexed
		\STATE \qquad \textbf{end}
		\STATE \qquad \texttt{\_syncthreads()}
		\STATE
		\STATE // Set $s$ to the largest number of power $2$ that is smaller than $n$
		\STATE \qquad s = \texttt{prevpow$(2,n-1)$}
		\STATE \qquad \textbf{while} $s > 1$
		\STATE \qquad \qquad \textbf{if} $(j \le s \ \&\& \ j+s \le n)$ \textbf{then}
		\STATE \qquad \qquad \qquad \texttt{shared\_mem[j]} $\leftarrow$ \texttt{shared\_mem[j]} $+$ \texttt{shared\_mem[j + s]}
		\STATE \qquad \qquad \textbf{end if}
		\STATE \qquad \qquad \texttt{\_syncthreads()}
		\STATE \qquad \qquad $s>>1$		\qquad //Divide $s$ by $2$
		\STATE \qquad \textbf{end for}
		
		\STATE \qquad \textbf{if} $j == 1$ \textbf{then}
		\STATE \qquad \qquad \texttt{norm2} $\leftarrow$ \texttt{shared\_mem[1]}
		\STATE \qquad \textbf{end if}
		\STATE \texttt{end}
	\end{algorithmic}\label{alg:norm2-parallel-reduction}
\end{algorithm}

\newpage
\section{Detailed benchmark results}\label{appendix:benchmark-table}

%\BenchmarkQPDetailTable{gpu_mittelmann_lp}{large LPs}
%\BenchmarkQPDetailTable{gpu_opf_large_lp}{large OPF LP}

\captionsetup{labelfont=bf}
\centering
\renewcommand{\detailtablecaption}{\bf Solve times and iteration counts for the large OPF SOCP problem set}
\begin{landscape}
\scriptsize
\begin{longtable}{||lcccc||cccc||cccc||}
\caption{\detailtablecaption}
\\
 & &  & & & \multicolumn{4}{c||}{\underline{iterations}}& \multicolumn{4}{c||}{\underline{total time (s)}}\\[2ex] 
Problem & vars. & cons. & nnz(A) & nnz(P)  & ClarabelGPU & Mosek* & ClarabelRs & Mosek & ClarabelGPU & Mosek* & ClarabelRs & Mosek\\[1ex]
\hline
\endhead
\sc{case2383wp\_k} & 97600 & 20393 & 155950 & 0 &   110 &   - &  \winner  107 &   - &  \winner  1.48 &   - &   2.79 &   -\\ 
\sc{case2383wp\_k\_\_api} & 97600 & 20393 & 155950 & 0 &  \winner  21 &   26 &   22 &   33 &  \winner  0.626 &   0.892 &   0.868 &   2.82\\ 
\sc{case2383wp\_k\_\_sad} & 97600 & 20393 & 155950 & 0 &  \winner  96 &   - &   97 &   - &  \winner  1.12 &   - &   2.94 &   -\\ 
\sc{case2312\_goc} & 98556 & 20518 & 159516 & 0 &  \winner  47 &   152 &  \winner  47 &   - &  \winner  0.802 &   4.02 &   1.58 &   -\\ 
\sc{case2312\_goc\_\_api} & 98556 & 20518 & 159516 & 0 &   65 &   244 &  \winner  64 &   - &  \winner  1.04 &   6.15 &   1.71 &   -\\ 
\sc{case2312\_goc\_\_sad} & 98556 & 20518 & 159516 & 0 &  \winner  59 &   191 &  \winner  59 &   - &  \winner  0.992 &   5.08 &   1.61 &   -\\ 
\sc{case2737sop\_k} & 109822 & 22777 & 176602 & 0 &  \winner  74 &   - &  \winner  74 &   - &  \winner  1.06 &   - &   3.66 &   -\\ 
\sc{case2737sop\_k\_\_api} & 109822 & 22777 & 176602 & 0 &  \winner  71 &   - &  \winner  71 &   - &  \winner  1.01 &   - &   2.37 &   -\\ 
\sc{case2737sop\_k\_\_sad} & 109822 & 22777 & 176602 & 0 &   72 &   - &  \winner  71 &   - &  \winner  1.03 &   - &   3.07 &   -\\ 
\sc{case2736sp\_k} & 110022 & 22878 & 176901 & 0 &  \winner  69 &   - &   74 &   - &  \winner  1.24 &   - &   2.15 &   -\\ 
\sc{case2736sp\_k\_\_api} & 110022 & 22878 & 176901 & 0 &  \winner  57 &   - &  \winner  57 &   - &  \winner  0.929 &   - &    1.9 &   -\\ 
\sc{case2736sp\_k\_\_sad} & 110022 & 22878 & 176901 & 0 &  \winner  72 &   - &  \winner  72 &   - &  \winner  1.01 &   - &    4.2 &   -\\ 
\sc{case2746wp\_k} & 111106 & 23320 & 178556 & 0 &   73 &   - &  \winner  71 &   - &  \winner  1.07 &   - &   2.47 &   -\\ 
\sc{case2746wp\_k\_\_api} & 111106 & 23320 & 178556 & 0 &  \winner  33 &   35 &   - &   34 &  \winner  0.76 &   1.26 &   - &   2.24\\ 
\sc{case2746wp\_k\_\_sad} & 111106 & 23320 & 178556 & 0 &  \winner  74 &   - &   75 &   - &  \winner  1.07 &   - &   2.26 &   -\\ 
\sc{case2746wop\_k} & 111822 & 23434 & 179829 & 0 &  \winner  65 &   - &  \winner  65 &   - &  \winner     1 &   - &   2.14 &   -\\ 
\sc{case2746wop\_k\_\_api} & 111822 & 23434 & 179829 & 0 &   28 &   31 &  \winner  25 &   36 &  \winner  0.913 &   1.16 &   1.22 &   2.44\\ 
\sc{case2746wop\_k\_\_sad} & 111822 & 23434 & 179829 & 0 &   64 &   - &  \winner  62 &   - &  \winner  0.966 &   - &   1.97 &   -\\ 
\sc{case3012wp\_k} & 120676 & 25202 & 193907 & 0 &  \winner  86 &   - &   89 &   - &  \winner  1.16 &   - &   5.32 &   -\\ 
\sc{case3012wp\_k\_\_api} & 120676 & 25202 & 193907 & 0 &  \winner  28 &   57 &   - &   34 &  \winner  0.77 &    2.1 &   - &   2.38\\ 
\sc{case3012wp\_k\_\_sad} & 120676 & 25202 & 193907 & 0 &   99 &   - &  \winner  94 &   - &  \winner  1.31 &   - &   2.86 &   -\\ 
\sc{case2848\_rte} & 122708 & 25858 & 197228 & 0 &  \winner  81 &   - &  \winner  81 &   - &  \winner   1.3 &   - &    3.4 &   -\\ 
\sc{case2848\_rte\_\_api} & 122708 & 25858 & 197228 & 0 &   76 &   218 &  \winner  75 &   - &  \winner  1.35 &    6.1 &    3.7 &   -\\ 
\sc{case2848\_rte\_\_sad} & 122708 & 25858 & 197228 & 0 &  \winner  70 &   217 &  \winner  70 &   - &  \winner  1.25 &      6 &   3.23 &   -\\ 
\sc{case2868\_rte} & 123912 & 26164 & 198869 & 0 &  \winner  84 &   228 &   86 &   - &  \winner  1.41 &   6.72 &   8.48 &   -\\ 
\sc{case2868\_rte\_\_api} & 123912 & 26164 & 198869 & 0 &  \winner  72 &   258 &   77 &   - &  \winner  1.55 &   7.62 &   3.84 &   -\\ 
\sc{case2868\_rte\_\_sad} & 123912 & 26164 & 198869 & 0 &  \winner  79 &   287 &  \winner  79 &   - &  \winner  1.32 &   8.18 &   6.56 &   -\\ 
\sc{case3120sp\_k} & 124354 & 25856 & 199827 & 0 &  \winner  82 &   - &   95 &   - &  \winner  1.21 &   - &   4.85 &   -\\ 
\sc{case3120sp\_k\_\_api} & 124354 & 25856 & 199827 & 0 &   96 &   - &  \winner  93 &   - &  \winner  1.32 &   - &   5.08 &   -\\ 
\sc{case2853\_sdet} & 128886 & 27445 & 206896 & 0 &  \winner  103 &   - &   105 &   - &  \winner  1.29 &   - &   6.31 &   -\\ 
\sc{case2853\_sdet\_\_api} & 128886 & 27445 & 206896 & 0 &   120 &   - &  \winner  116 &   - &  \winner  1.41 &   - &   3.45 &   -\\ 
\sc{case2853\_sdet\_\_sad} & 128886 & 27445 & 206896 & 0 &  \winner  102 &   278 &  \winner  102 &   - &  \winner  1.27 &   8.68 &   3.01 &   -\\ 
\sc{case3022\_goc} & 134780 & 28060 & 217434 & 0 &  \winner  41 &   176 &   - &   140 &  \winner  0.984 &   5.86 &   - &   9.95\\ 
\sc{case3022\_goc\_\_api} & 134780 & 28060 & 217434 & 0 &   55 &   203 &  \winner  54 &   170 &  \winner   1.1 &    6.8 &    2.2 &   12.3\\ 
\sc{case3022\_goc\_\_sad} & 134780 & 28060 & 217434 & 0 &  \winner  41 &   177 &  \winner  41 &   141 &  \winner  0.971 &   5.76 &   1.74 &   9.79\\ 
\sc{case3375wp\_k} & 139126 & 29112 & 223999 & 0 &   86 &   - &  \winner  85 &   - &  \winner  1.26 &   - &   3.05 &   -\\ 
\sc{case3375wp\_k\_\_api} & 139126 & 29112 & 223999 & 0 &   161 &   - &  \winner  147 &   - &  \winner  1.86 &   - &   4.82 &   -\\ 
\sc{case3375wp\_k\_\_sad} & 139126 & 29112 & 223999 & 0 &   104 &   - &  \winner  93 &   - &  \winner   1.4 &   - &   3.48 &   -\\ 
\sc{case2869\_pegase} & 143608 & 30153 & 236217 & 0 &  \winner  53 &   236 &  \winner  53 &   - &  \winner  1.25 &   8.46 &    4.9 &   -\\ 
\sc{case2869\_pegase\_\_api} & 143608 & 30153 & 236217 & 0 &  \winner  50 &   262 &  \winner  50 &   - &  \winner  1.24 &   9.31 &   6.44 &   -\\ 
\sc{case2869\_pegase\_\_sad} & 143608 & 30153 & 236217 & 0 &  \winner  48 &   278 &  \winner  48 &   - &  \winner  1.22 &   9.66 &   2.99 &   -\\ 
\sc{case2742\_goc} & 144110 & 29856 & 236467 & 0 &  \winner  109 &   118 &   113 &   - &  \winner  2.31 &   4.81 &   14.4 &   -\\ 
\sc{case2742\_goc\_\_api} & 144110 & 29856 & 236467 & 0 &  \winner  101 &   128 &   103 &   - &  \winner  2.15 &   5.06 &    7.2 &   -\\ 
\sc{case2742\_goc\_\_sad} & 144110 & 29856 & 236467 & 0 &  \winner  106 &   118 &   108 &   - &  \winner  2.32 &   4.67 &   7.47 &   -\\ 
\sc{case4661\_sdet} & 198498 & 41599 & 320728 & 0 &  \winner  139 &   - &   140 &   - &  \winner   2.7 &   - &   16.8 &   -\\ 
\sc{case4661\_sdet\_\_api} & 198498 & 41599 & 320728 & 0 &   126 &   - &  \winner  124 &   - &  \winner  2.59 &   - &   10.7 &   -\\ 
\sc{case4661\_sdet\_\_sad} & 198498 & 41599 & 320728 & 0 &   140 &   322 &  \winner  138 &   - &  \winner  2.71 &   18.6 &    8.3 &   -\\ 
\sc{case3970\_goc} & 205819 & 42789 & 337328 & 0 &   140 &   156 &  \winner  139 &   - &  \winner  3.79 &   9.51 &     12 &   -\\ 
\sc{case3970\_goc\_\_api} & 205819 & 42789 & 337328 & 0 &  \winner  114 &   137 &   116 &   - &  \winner  3.59 &   8.41 &   12.3 &   -\\ 
\sc{case3970\_goc\_\_sad} & 205819 & 42789 & 337328 & 0 &   147 &   168 &  \winner  141 &   - &  \winner  3.23 &   10.2 &     12 &   -\\ 
\sc{case4020\_goc} & 216455 & 44877 & 355706 & 0 &   187 &   - &  \winner  180 &   - &  \winner  4.64 &   - &   18.3 &   -\\ 
\sc{case4020\_goc\_\_api} & 216455 & 44877 & 355706 & 0 &  \winner  98 &   142 &   101 &   - &  \winner  3.13 &   8.86 &   21.7 &   -\\ 
\sc{case4020\_goc\_\_sad} & 216455 & 44877 & 355706 & 0 &  \winner  180 &   - &   187 &   - &  \winner  4.57 &   - &   20.6 &   -\\ 
\sc{case4917\_goc} & 218213 & 45522 & 352833 & 0 &  \winner  52 &   166 &  \winner  52 &   147 &  \winner   1.7 &   8.82 &   4.82 &   16.5\\ 
\sc{case4917\_goc\_\_api} & 218213 & 45522 & 352833 & 0 &  \winner  62 &   256 &  \winner  62 &   - &  \winner  1.86 &   13.4 &   5.98 &   -\\ 
\sc{case4917\_goc\_\_sad} & 218213 & 45522 & 352833 & 0 &  \winner  48 &   183 &   - &   169 &  \winner  1.69 &    9.8 &   - &   18.7\\ 
\sc{case4601\_goc} & 225588 & 46885 & 368445 & 0 &  \winner  180 &   192 &  \winner  180 &   - &  \winner  3.84 &   11.7 &   19.8 &   -\\ 
\sc{case4601\_goc\_\_api} & 225588 & 46885 & 368445 & 0 &  \winner  109 &   148 &  \winner  109 &   156 &  \winner  3.03 &   8.96 &   14.4 &   19.3\\ 
\sc{case4601\_goc\_\_sad} & 225588 & 46885 & 368445 & 0 &   179 &   191 &  \winner  178 &   - &  \winner  4.56 &   12.2 &     15 &   -\\ 
\sc{case4837\_goc} & 240160 & 49855 & 392884 & 0 &  \winner  134 &   154 &   138 &   - &  \winner  3.82 &   10.3 &   36.1 &   -\\ 
\sc{case4837\_goc\_\_api} & 240160 & 49855 & 392884 & 0 &  \winner  123 &   177 &  \winner  123 &   - &  \winner  3.73 &   12.1 &   13.7 &   -\\ 
\sc{case4837\_goc\_\_sad} & 240160 & 49855 & 392884 & 0 &   133 &   162 &  \winner  131 &   - &  \winner  3.88 &     11 &   17.5 &   -\\ 
\sc{case4619\_goc} & 254753 & 52616 & 419210 & 0 &   147 &   151 &  \winner  145 &   160 &  \winner  4.51 &   10.7 &   19.9 &   22.2\\ 
\sc{case4619\_goc\_\_api} & 254753 & 52616 & 419210 & 0 &  \winner  129 &   161 &   131 &   - &  \winner  4.01 &   11.7 &   15.8 &   -\\ 
\sc{case4619\_goc\_\_sad} & 254753 & 52616 & 419210 & 0 &  \winner  148 &   151 &   150 &   - &  \winner  5.15 &     11 &   19.2 &   -\\ 
\sc{case5658\_epigrids} & 282180 & 58620 & 461724 & 0 &   146 &   177 &  \winner  143 &   - &  \winner  5.53 &   12.9 &   19.7 &   -\\ 
\sc{case5658\_epigrids\_\_api} & 282180 & 58620 & 461724 & 0 &   163 &   203 &  \winner  162 &   - &  \winner  5.45 &   14.6 &   25.3 &   -\\ 
\sc{case5658\_epigrids\_\_sad} & 282180 & 58620 & 461724 & 0 &   163 &   217 &  \winner  157 &   - &  \winner   4.8 &   15.7 &   20.4 &   -\\ 
\sc{case6468\_rte} & 286248 & 59396 & 463599 & 0 &   103 &   250 &  \winner  98 &   - &  \winner  3.01 &   15.8 &   12.4 &   -\\ 
\sc{case6468\_rte\_\_api} & 286248 & 59396 & 463599 & 0 &  \winner  88 &   - &   90 &   - &  \winner  3.59 &   - &   15.2 &   -\\ 
\sc{case6468\_rte\_\_sad} & 286248 & 59396 & 463599 & 0 &  \winner  98 &   255 &   99 &   - &  \winner  3.08 &   16.4 &   19.7 &   -\\ 
\sc{case6470\_rte} & 287806 & 60144 & 465757 & 0 &   113 &   321 &  \winner  111 &   - &  \winner  3.29 &   20.5 &   26.4 &   -\\ 
\sc{case6470\_rte\_\_api} & 287806 & 60144 & 465757 & 0 &  \winner  85 &   327 &  \winner  85 &   - &  \winner  2.95 &   21.1 &   15.5 &   -\\ 
\sc{case6470\_rte\_\_sad} & 287806 & 60144 & 465757 & 0 &   104 &   354 &  \winner  100 &   - &  \winner  3.66 &   22.2 &   14.4 &   -\\ 
\sc{case6495\_rte} & 288050 & 60099 & 465939 & 0 &  \winner  114 &   323 &  \winner  114 &   - &  \winner  3.37 &   20.4 &   16.4 &   -\\ 
\sc{case6495\_rte\_\_api} & 288050 & 60099 & 465939 & 0 &   - &   343 &  \winner  98 &   - &   - &   21.4 &  \winner  12.4 &   -\\ 
\sc{case6495\_rte\_\_sad} & 288050 & 60099 & 465939 & 0 &   112 &   344 &  \winner  108 &   - &  \winner  3.25 &   22.1 &   25.2 &   -\\ 
\sc{case6515\_rte} & 288710 & 60239 & 466922 & 0 &   97 &   294 &  \winner  95 &   - &  \winner  3.19 &   18.3 &   12.6 &   -\\ 
\sc{case6515\_rte\_\_api} & 288710 & 60239 & 466922 & 0 &   85 &   341 &  \winner  81 &   - &  \winner  2.97 &   21.4 &   21.4 &   -\\ 
\sc{case6515\_rte\_\_sad} & 288710 & 60239 & 466922 & 0 &  \winner  84 &   324 &   88 &   - &  \winner  2.79 &   20.4 &   12.1 &   -\\ 
\sc{case7336\_epigrids} & 358450 & 74618 & 585766 & 0 &   155 &   195 &  \winner  153 &   - &  \winner  6.51 &   17.2 &   35.9 &   -\\ 
\sc{case7336\_epigrids\_\_api} & 358450 & 74618 & 585766 & 0 &   152 &   206 &  \winner  151 &   - &  \winner   6.1 &   18.4 &   32.8 &   -\\ 
\sc{case7336\_epigrids\_\_sad} & 358450 & 74618 & 585766 & 0 &  \winner  156 &   188 &  \winner  156 &   - &  \winner  6.56 &   17.8 &   58.6 &   -\\ 
\sc{case10000\_goc} & 440997 & 92799 & 712177 & 0 &  \winner  71 &   87 &   - &   - &  \winner  3.58 &   9.49 &   - &   -\\ 
\sc{case10000\_goc\_\_api} & 440997 & 92799 & 712177 & 0 &  \winner  62 &   205 &   63 &   - &  \winner  3.97 &   20.4 &   15.6 &   -\\ 
\sc{case10000\_goc\_\_sad} & 440997 & 92799 & 712177 & 0 &  \winner  89 &   152 &   - &   - &  \winner  4.13 &     16 &   - &   -\\ 
\sc{case8387\_pegase} & 459046 & 96351 & 750588 & 0 &  \winner  137 &   - &   138 &   - &  \winner   4.8 &   - &   22.3 &   -\\ 
\sc{case8387\_pegase\_\_api} & 459046 & 96351 & 750588 & 0 &  \winner  141 &   - &  \winner  141 &   - &  \winner  5.08 &   - &   25.9 &   -\\ 
\sc{case8387\_pegase\_\_sad} & 459046 & 96351 & 750588 & 0 &   - &   - &  \winner  142 &   - &   - &   - &  \winner  34.9 &   -\\ 
\sc{case9591\_goc} & 495073 & 102120 & 811889 & 0 &   229 &  \winner  206 &   226 &   - &  \winner    12 &   32.7 &   69.9 &   -\\ 
\sc{case9591\_goc\_\_api} & 495073 & 102120 & 811889 & 0 &  \winner  149 &   162 &  \winner  149 &   190 &  \winner  7.92 &   26.3 &   41.6 &   49.3\\ 
\sc{case9591\_goc\_\_sad} & 495073 & 102120 & 811889 & 0 &   230 &   - &  \winner  227 &   - &  \winner  10.8 &   - &   78.6 &   -\\ 
\sc{case9241\_pegase} & 502110 & 104741 & 827419 & 0 &   72 &   - &  \winner  67 &   - &  \winner  5.29 &   - &   14.9 &   -\\ 
\sc{case9241\_pegase\_\_api} & 502110 & 104741 & 827419 & 0 &  \winner  69 &   - &   - &   - &  \winner  4.02 &   - &   - &   -\\ 
\sc{case9241\_pegase\_\_sad} & 502110 & 104741 & 827419 & 0 &   72 &   - &  \winner  68 &   - &  \winner  4.26 &   - &   21.7 &   -\\ 
\sc{case10192\_epigrids} & 529695 & 109758 & 866501 & 0 &  \winner  94 &   113 &   95 &   - &  \winner  5.74 &   16.5 &   29.5 &   -\\ 
\sc{case10192\_epigrids\_\_api} & 529695 & 109758 & 866501 & 0 &  \winner  93 &   115 &   - &   - &  \winner  5.92 &   16.7 &   - &   -\\ 
\sc{case10192\_epigrids\_\_sad} & 529695 & 109758 & 866501 & 0 &  \winner  110 &   - &   114 &   - &  \winner  6.36 &   - &   33.2 &   -\\ 
\sc{case10480\_goc} & 573754 & 118760 & 943278 & 0 &  \winner  160 &   235 &   161 &   - &  \winner  9.62 &   42.4 &   53.5 &   -\\ 
\sc{case10480\_goc\_\_api} & 573754 & 118760 & 943278 & 0 &  \winner  139 &   213 &   - &   - &  \winner  8.91 &   39.3 &   - &   -\\ 
\sc{case10480\_goc\_\_sad} & 573754 & 118760 & 943278 & 0 &  \winner  162 &   - &   164 &   - &  \winner    11 &   - &   54.3 &   -\\ 
\sc{case13659\_pegase} & 662910 & 140961 & 1081042 & 0 &  \winner  80 &   - &   82 &   - &  \winner  6.53 &   - &   30.8 &   -\\ 
\sc{case13659\_pegase\_\_api} & 662910 & 140961 & 1081042 & 0 &  \winner  78 &   - &  \winner  78 &   - &  \winner  6.35 &   - &   29.2 &   -\\ 
\sc{case13659\_pegase\_\_sad} & 662910 & 140961 & 1081042 & 0 &   78 &   331 &  \winner  77 &   - &  \winner  6.29 &   59.7 &   30.5 &   -\\ 
\sc{case20758\_epigrids} & 1048059 & 218151 & 1709288 & 0 &  \winner  123 &   - &   - &   - &  \winner  12.4 &   - &   - &   -\\ 
\sc{case20758\_epigrids\_\_api} & 1048059 & 218151 & 1709288 & 0 &   153 &   - &  \winner  150 &   - &  \winner    16 &   - &    112 &   -\\ 
\sc{case20758\_epigrids\_\_sad} & 1048059 & 218151 & 1709288 & 0 &  \winner  136 &   - &  \winner  136 &   - &  \winner  15.6 &   - &   79.9 &   -\\ 
\sc{case19402\_goc} & 1064421 & 219911 & 1753874 & 0 &   259 &   - &  \winner  254 &   - &  \winner  21.8 &   - &    168 &   -\\ 
\sc{case19402\_goc\_\_api} & 1064421 & 219911 & 1753874 & 0 &  \winner  215 &   - &  \winner  215 &   - &  \winner  19.3 &   - &    165 &   -\\ 
\sc{case19402\_goc\_\_sad} & 1064421 & 219911 & 1753874 & 0 &  \winner  253 &   - &  \winner  253 &   - &  \winner  21.5 &   - &    175 &   -\\ 
\sc{case30000\_goc\_\_api} & 1195834 & 249462 & 1921390 & 0 &  \winner  146 &   - &   - &   - &  \winner  13.1 &   - &   - &   -\\ 
\sc{case24464\_goc} & 1202964 & 248644 & 1983519 & 0 &   233 &   - &  \winner  219 &   - &  \winner  21.2 &   - &    133 &   -\\ 
\sc{case24464\_goc\_\_api} & 1202964 & 248644 & 1983519 & 0 &   227 &   - &  \winner  210 &   - &  \winner  16.9 &   - &    108 &   -\\ 
\sc{case24464\_goc\_\_sad} & 1202964 & 248644 & 1983519 & 0 &   226 &   - &  \winner  221 &   - &  \winner  17.9 &   - &    124 &   -\\ 
\sc{case78484\_epigrids} & 3904758 & 811998 & 6389202 & 0 &   413 &   - &  \winner  410 &   - &  \winner   118 &   - &   1.02e+03 &   -\\ 
\sc{case78484\_epigrids\_\_api} & 3904758 & 811998 & 6389202 & 0 &  \winner  396 &   - &   - &   - &  \winner   138 &   - &   - &   -\\ 
\sc{case78484\_epigrids\_\_sad} & 3904758 & 811998 & 6389202 & 0 &  \winner  405 &   - &  \winner  405 &   - &  \winner   103 &   - &    981 &   -\\ 
\end{longtable}

\end{landscape}
%\label{table:#1:detail}

\end{document}